\documentclass[twoside]{amsart}
\usepackage[american]{babel}
\usepackage{amsopn,amssymb}
\usepackage{fancyhdr}
\theoremstyle{plain}
\newtheorem{thm}{Theorem}[section]
\newtheorem{lmm}[thm]{Lemma}
\newtheorem{crllr}[thm]{Corollary}
\newtheorem{prpstn}[thm]{Proposition}

\newtheorem{rmrk}[thm]{Remark}
\def\xR{{\mathbb R}}

\def\xN{{\mathbb N}}

\def\p{{\partial}}
\def\xL{{\rm L}}
\def\xW{{\rm W}}
\def\xCn{{\rm C}}
\def\xdif{{\rm d}}
\def\xLtwo{{{\rm L}^2}}
\def\xHone{{{\rm H}^1}}

\def\xLn{{\rm L}}
\def\xHn{{\rm H}}

\def\xim{{\rm im}}
\def\xker{{\rm ker}}

\def\pim{{\mathfrak{Im}}}
\def\pre{{\mathfrak{Re}}}

\title{Exit from a basin of attraction for stochastic weakly damped nonlinear Schr\"odinger equations}

\author{Eric GAUTIER$^{1,2}$}
\begin{document}
\maketitle
\begin{abstract}
\ We consider weakly damped nonlinear Schr\"odinger equations
perturbed by a noise of small amplitude.\@ The small noise is either
complex and of additive type or real and of multiplicative type.\@
It is white in time and colored in space.\@ Zero is an
asymptotically stable equilibrium point of the deterministic
equations.\@ We study the exit from a neighborhood of zero,
invariant by the flow of the deterministic equation, in $\xLtwo$ or
in $\xHone$.\@ Due to noise, large fluctuations off zero occur.\@
Thus, on a sufficiently large time scale, exit from these domains of
attraction occur.\@ A formal characterization of the small noise
asymptotic of both the first exit times and the exit points is
given.\vspace{0.2cm}

\vspace{0.2cm}

 \vspace{0.2cm}
\end{abstract}

\footnotetext[1]{IRMAR, Ecole Normale Sup\'erieure de Cachan,
antenne de Bretagne, Campus de Ker Lann, avenue R. Schuman, 35170
Bruz, France} \footnotetext[2]{CREST-INSEE, URA D2200, 3 avenue
Pierre Larousse, 92240 Malakoff, France}
\pagestyle{myheadings}\markboth{{A. Debussche and E.
Gautier}}{{Small noise asymptotic of the timing jitter in soliton
transmission}}

\section{Introduction}\label{s1}
The study of the first exit time from a neighborhood of an
asymptotically stable equilibrium point, the exit place
determination or the transition between two equilibrium points in
randomly perturbed dynamical systems is important in several areas
of physics among which statistical and quantum mechanics, chemical
reactions, the natural sciences, macroeconomics as to model currency
crises or escape in learning models...\\
\indent For a fixed noise amplitude and for diffusions, the first
exit time and the distribution of the exit points on the boundary of
a domain can be characterized respectively by the Dirichlet and
Poisson equations.\@ However, when the dimension is larger than one,
we may seldom solve explicitly these equations and large deviation
techniques are precious tools when the noise is assumed to be small;
see for example \cite{DZ,FW}.\@ The techniques used in the physics
literature is often called optimal fluctuations or instanton
formalism and are closely related to large deviations.\\
\indent In that case, an energy generally characterizes the
transition between two states and the exit from a neighborhood of an
asymptotically stable equilibrium point of the deterministic
equation.\@ The energy is derived from the rate function of the
sample path large deviation principle (LDP).\@ When a LDP holds, the
first order of the probability of rare events is that of the
Boltzman theory and the square of the amplitude of the small noise
acts as the temperature.\@ The deterministic dynamics is sometimes
interpreted as the evolution at temperature 0 and the small noise as
the small temperature nonequilibrium case.\@ The exit or transition
problem is then related to a deterministic least-action principle.\@
The paths that minimize the energy, also called minimum action
paths, are the most likely exiting paths or transitions.\@ When the
infimum is unique, the system has a behavior which is almost
deterministic even though there is noise.\@ Indeed, other possible
exiting paths, points or transitions are exponentially less
probable.\@ In the pioneering article \cite{FJL}, a nonlinear heat
equation perturbed by a small noise of additive type is
considered.\@ Transitions in that case prove to be the instantons of
quantum mechanics.\@ The problem is studied again in \cite{ERVE}
where a numerical scheme is presented to compute the optimal
paths.\@ In \cite{LMMCS}, mathematical and numerical
predictions for a noisy exit problem are confirmed experimentally.\\
\indent In this article, we consider the case of weakly damped
nonlinear Schr\"odinger (NLS) equations in $\xR^d$.\@ These
equations are a generic model for the propagation of the enveloppe
of a wave packet in weakly nonlinear and dispersive media.\@ They
appear for example in nonlinear optics, hydrodynamics, biology,
field theory, crystals Fermi-Pasta-Ulam chains of atoms.\@ The
equations are perturbed by a small noise.\@ In optics, the noise
corresponds to the spontaneous emission noise due to amplifiers
placed along the fiber line in order to compensate for loss,
corresponding to the weak damping, in the fiber.\@ We shall consider
here that there remains a small weak damping term.\@ In the context
of crystals or of Fermi-Pasta-Ulam chains of atoms, the noise
accounts for thermal effects.\@ The relevance of the study of the
exit from a domain in nonlinear optics is discussed in
\cite{KLMMC}.\@ The noise is of additive or multiplicative type.\@
We define it as the time derivative in the sense of distributions of
a Hilbert space-valued Wiener process $\left(W_t\right)_{t\geq0}$.\@
The evolution equation could be written in It\^o form
\begin{equation}\label{e1}
i \xdif u^{\epsilon,u_0}=(\Delta u^{\epsilon,u_0} + \lambda
|u^{\epsilon,u_0}|^{2\sigma}u^{\epsilon,u_0}-i\alpha
u^{\epsilon,u_0})\xdif t + \sqrt{\epsilon}\xdif W,
\end{equation}
where $\alpha$ and $\epsilon$ are positive and $u_0$ is an initial
datum in $\xLtwo$ or $\xHone$.\@ When the noise is of multiplicative
type, the product is a Stratonovich product and the equation may be
written
\begin{equation}\label{e1m}
i \xdif u^{\epsilon,u_0}=(\Delta u^{\epsilon,u_0} + \lambda
|u^{\epsilon,u_0}|^{2\sigma}u^{\epsilon,u_0}-i\alpha
u^{\epsilon,u_0})\xdif t + \sqrt{\epsilon}u^{\epsilon,u_0}\circ\xdif
W.
\end{equation}
Contrary to the Heat equation the linear part has no smoothing
effects.\@ In our case, it defines a linear group which is an
isometry on the $\xLtwo$ based Sobolev spaces.\@ Thus, we cannot
treat spatially rough noises and consider colored in space Wiener
processes.\@ This latter property is required to obtain {\it
bona-fide} Wiener processes in infinite dimensions.\@ The white
noise often considered in Physics seems to give rise to ill-posed
problems.\\
\indent Results on local and global well-posedness and on the effect
of a noise on the blow-up phenomenon are proved in
\cite{dBD0,dBD1,dBD2,dBD3} in the case $\alpha=0$.\@ Mixing property
and convergence to equilibrium is studied for weakly damped cubic
one dimensional equations on a bounded domain in \cite{DO}.\@ We
consider these equations in the whole space $\mathbb{R}^d$ and
assume that the power of the nonlinearity $\sigma$ satisfies
$\sigma<2/d$. We may check that the above result still hold with the
damping term and that for such powers of the nonlinearity the
solutions
do not exhibit blow-up.\\
\indent In \cite{EG1} and \cite{EG2}, we have proved sample paths
LDPs for the two types of noises but without damping and deduced the
asymptotic of the tails of the blow-up times.\@ In \cite{EG1}, we
also deduced the tails of the mass, defined later, of the pulse at
the end of a fiber optical line.\@ We have thus evaluated the error
probabilities in optical soliton transmission when the receiver
records the signal on an infinite time interval.\@ In \cite{DG} we
have applied the LDPs to the problem of the diffusion in position of
the soliton and studied the tails of the random arrival time of a
pulse in optical soliton transmission for noises of additive and
multiplicative types.\\
\indent The flow defined by the above equations can be decomposed in
a Hamiltonian, a gradient and a random component.\@ The mass
\begin{equation*}
\mathbf{N}\left(u\right)=\int_{\mathbb{R}^d}\left| u\right|^2 dx
\end{equation*}
characterizes the gradient component.\@ The Hamiltonian denoted by
$\mathbf{H}(u)$, defined for functions in $\xHone$, has a kinetic
and a potential term, it may be written
\begin{equation*}
\mathbf{H}\left(u\right)=(1/2)\int_{\mathbb{R}^d}\left|\nabla
u\right|^2
dx-(\lambda/(2\sigma+2))\int_{\mathbb{R}^d}\left|u\right|^{2\sigma+2}dx.
\end{equation*}
Note that the vector fields associated to the mass and Hamiltonian
are orthogonal.\@ We could rewrite, for example equation \eqref{e1},
as
\begin{equation*}
\xdif u^{\epsilon,u_0}=\left(\frac{\delta \mathbf{H}
\left(u^{\epsilon,u_0}\right)}{\delta \overline{u^{\epsilon,u_0}}}-
(\alpha/2)\frac{\delta
\mathbf{N}\left(u^{\epsilon,u_0}\right)}{\delta
u^{\epsilon,u_0}}\right)\xdif t -i \sqrt{\epsilon}\xdif W.
\end{equation*}
Also, the mass and Hamiltonian are invariant quantities of the
equation without noise and damping.\@ Other quantities like the
linear or angular momentum are also invariant for nonlinear
Schr\"odinger equations.\\
\indent Without noise, solutions are uniformly attracted to zero in
$\xLtwo$ and in $\xHone$.\@ In this article we study the classical
problem of exit from a bounded domain containing zero in its
interior and invariant by the deterministic evolution.\@ We prove
that the behavior of the random evolution is completely different
from the deterministic evolution.\@ Though for finite times the
probabilities of large excursions off neighborhoods of zero go to
zero exponentially fast with $\epsilon$, if we wait long enough -
the time scale is exponential - such large fluctuations occur and
exit from a domain takes place.\@ We give two types of results
depending on the topology we consider, $\xLtwo$ or $\xHone$.\@ The
$\xLtwo$-setting is less involved than the $\xHone$-setting.\@ This
is due to the structure of the NLS equation and the fact that the
$\xLtwo$ norm is conserved for deterministic non damped equations.\@
We have chosen to also work in $\xHone$ because it is the
mathematical framework to study perturbations of solitons; a problem
we hope to address in future research.\\
\indent We give a formal characterization of the small noise
asymptotic of the first exit time and exit points.\@ The main tool
is a uniform large deviation principle at the level of the paths of
the solutions.\@ The behavior of the process is proved to be
exponentially equivalent to that of the process starting from a
little ball around zero.\@ Thus, if a multiplicative noise and the
$\xLtwo$ topology is considered such balls are invariant by the
stochastic evolution as well and the exit problem is not
interesting.\@ In infinite dimensions we are faced with two major
difficulties.\@ Primarily, the domains under consideration are not
relatively compact.\@ In bounded domains of $\xR^d$, it is sometimes
possible to use compact embedding and the regularizing properties of
the semi-group.\@ In \cite{Fr} where the case of the Heat semi-group
and a space variable in a unidimensional torus is treated, these
properties are at hand.\@ Also, in \cite{CM}, the neighborhood is
defined for a strong topology of $\beta$-H\"older functions and is
relatively compact for a weaker topology, the space variable is
again in a bounded subset of $\xR^d$.\@ We are not able to use the
above properties here since the Schr\"odinger linear group is an
isometry on every Sobolev space based on $\xLtwo$ and we work on the
whole space $\xR^d$.\@ Another difficulty in infinite dimensions and
with unbounded linear operators is that, unlike ODEs, continuity of
the linear flow with respect to the initial data holds in a weak
sense.\@ The semi-group is strongly continuous and not in general
uniformly continuous.\@ We see that we may use other arguments than
those used in the finite dimensional setting, some of which are
taken from \cite{DPZ}, and that the expected results still hold.\@
We are also faced with particular difficulties arising from the
nonlinear Schr\"odinger equation among which the fact that the
nonlinearity is locally Lipschitz only in $\xHone$ for $d=1$.\@ In
this purpose, we use the hyper contractivity governed by the
Strichartz inequalities which is related to the dispersive
properties of the equation.\\
\indent In this article, we do not address the control problems for
the controlled deterministic PDE.\@ We could expect that the upper
and lower bound on the expected first exit time are equal and could
be written in terms of the usual quasi-potential. The exit points could
be related to solitary waves. These issues will be studied in future works.\\
\indent The article is organized as follows.\@ In the first section,
we introduce the main notations and tools, the proof of the uniform
large deviation principle is given in the annex.\@ In the next
section, we consider the exit off a domain in $\xLtwo$ for equations
with additive noise while in the last section we consider the exit
off domains in $\xHone$ for equations with an additive or
multiplicative noise.\@

\section{Preliminaries}\label{s2}
Throughout the paper the following notations are used.\\
\indent The set of positive integers and positive real numbers are
denoted by $\xN^*$ and $\xR_+^*$.\@ For $p\in \xN^*$, $\xLn^{p}$ is
the Lebesgue space of complex valued functions.\@ For $k$ in
$\xN^*$, $\xW^{k,p}$ is the Sobolev space of $\xLn^{p}$ functions
with partial derivatives up to level $k$, in the sense of
distributions, in $\xLn^{p}$.\@ For $p=2$ and $s$ in $\xR_+^*$,
$\xHn^{s}$ is the Sobolev space of tempered distributions $v$ of
Fourier transform $\hat{v}$ such that $(1+|\xi|^2)^{s/2}\hat{v}$
belongs to $\xLtwo$.\@ We denote the spaces by $\xLn_{\xR}^p$,
$\xW_{\xR}^{k,p}$ and $\xHn_{\xR}^s$ when the functions are
real-valued.\@ The space $\xLtwo$ is endowed with the inner product
$(u,v)_{\xL2}=\pre\int_{\xR^d}u(x)\overline{v}(x)dx$.\@ If $I$ is an
interval of $\xR$, $(E,\|\cdot\|_E)$ a Banach space and $r$ belongs
to $[1,\infty]$, then $\xLn^{r}(I;E)$ is the space of strongly
Lebesgue measurable functions $f$ from $I$
into $E$ such that $t\rightarrow \|f(t)\|_E$ is in $\xLn^{r}(I)$.\\
\indent The space of linear continuous operators from $B$ into
$\tilde{B}$, where $B$ and $\tilde{B}$ are Banach spaces is
$\mathcal{L}_c\left(B,\tilde{B}\right)$.\@ When $B=H$ and
$\tilde{B}=\tilde{H}$ are Hilbert spaces, such an operator is
Hilbert-Schmidt when $\sum_{j\in\xN}\|\Phi
e^H_j\|_{\tilde{H}}^2<\infty$ for every $\left(e_j\right)_{j\in\xN}$
complete orthonormal system of $H$.\@ The set of such operators is
denoted by $\mathcal{L}_2(H,\tilde{H})$, or $\mathcal{L}_2^{s,r}$
when $H=\xHn^{s}$ and $\tilde{H}=\xHn^{r}$.\@ When
$H=\xHn_{\xR}^{s}$ and
 $\tilde{H}=\xHn_{\xR}^{r}$, we denote it by $\mathcal{L}_{2,\xR}^{s,r}$.\@
When $s=0$ or $r=0$ the Hilbert space is $\xLtwo$ or $\xLn_{\xR}^2$.\\
\indent We also denote by $B_{\rho}^0$ and $S_{\rho}^0$ respectively
the open ball and the sphere centered at $0$ of radius $\rho$ in
$\xLtwo$.\@ We denote these by $B_{\rho}^1$ and $S_{\rho}^1$ in
$\xHone$.\@ We write $\mathcal{N}^0\left(A,\rho\right)$ for the
$\rho-$neighborhood of a set $A$ in $\xLtwo$ and
$\mathcal{N}^1\left(A,\rho\right)$ the neighborhood in $\xHone$.\@
In the following we impose that compact
sets satisfy the Hausdorff property.\\
\indent We use in Lemma \ref{l1} below the integrability of the
Schr\"odinger linear group which is related to the dispersive
property.\@ Recall that $(r(p),p)$ is an admissible pair if $p$ is
such that $2\leq p<2d/(d-2)$ when $d>2$ ($2\leq p<\infty$ when $d=2$
and $2\leq p\leq\infty$ when $d=1$) and $r(p)$ satisfies
$2/r(p)=d\left(1/2-1/p\right)$.\\
For every $(r(p),p)$ admissible pair and $T$ positive, we define the
Banach spaces
\begin{equation*}
Y^{(T,p)}=\xCn\left([0,T];\xLtwo\right)\cap
\xLn^{r(p)}\left(0,T;\xLn^{p}\right),
\end{equation*}
and
\begin{equation*}
X^{(T,p)}=\xCn\left([0,T];\xHone\right)\cap
\xLn^{r(p)}\left(0,T;\xW^{1,p}\right),
\end{equation*}
where the norms are the maximum of the norms in the two intersected
Banach spaces.\@ The Schr\"odinger linear group is denoted by
$\left(U(t)\right)_{t\geq0}$; it is defined on $\xLtwo$ or on
$\xHone$.\@ Let us recall the Strichartz inequalities, see
\cite{CAZ},\vspace{0.2cm}

\noindent\begin{tabular}{ll} (i)&There exists $C$ positive such that
for $u_0$ in $\xLtwo$, $T$ positive and\\& $(r(p),p)$
admissible pair,\\
&\\
&\hspace{2cm}$\left\|U(t)u_0\right\|_{Y^{(T,p)}}\leq C\left\|u_0\right\|_{\xLtwo},$\\
&\\
(ii)&For every $T$ positive, $(r(p),p)$ and $(r(q),q)$ admissible
pairs, $s$ and $\rho$\\
& such that $1/s+1/r(q)=1$ and
$1/\rho+1/q=1$, there exists $C$ positive such\\
&that for $f$ in $\xLn^s\left(0,T;\xLn^{\rho}\right)$,\\
&\\
&\hspace{2cm}$\left\|\int_0^{\cdot}U(\cdot-s)f(s)ds\right\|_{Y^{(T,p)}}\leq
C\|f\|_{\xLn^{s}\left(0,T;\xLn^{\rho}\right)}.$
\end{tabular}\vspace{0.3cm}

\noindent Similar inequalities hold when the group is acting on
$\xHone$, replacing $\xLtwo $ by $\xHone$, $Y^{(T,p)}$ by
$X^{(T,p)}$ and $\xLn^{s}\left(0,T;\xLn^{\rho}\right)$ by
$\xLn^{s}\left(0,T;\xW^{1,\rho}\right)$.\\ \indent It is known that,
in the Hilbert space setting, only direct images of uncorrelated
space wise Wiener processes by Hilbert-Schmidt operators are well
defined.\@ However, when the semi-group has regularizing properties,
the semi-group may act as a Hilbert-Schmidt operator and a white in
space noise may be considered.\@ It is not possible here since the
Schr\"odinger group is an isometry on the Sobolev spaces based on
$\xLtwo$.\@ The Wiener process $W$ is thus defined as $\Phi W_c$,
where $W_c$ is a cylindrical Wiener process on $\xLtwo$ and $\Phi$
is Hilbert-Schmidt.\@
Then $\Phi\Phi^*$ is the correlation operator of $W(1)$, it has finite trace.\\
\indent We consider the following Cauchy problems
\begin{equation}\label{ec1}
\left\{\begin{array}{rl} i \xdif u^{\epsilon,u_0}&=(\Delta
u^{\epsilon,u_0} + \lambda
|u^{\epsilon,u_0}|^{2\sigma}u^{\epsilon,u_0}-i\alpha
u^{\epsilon,u_0})\xdif t + \sqrt{\epsilon}\xdif W,\\
u^{\epsilon,u_0}(0)&=u_0
\end{array}\right.
\end{equation}
with $u_0$ in $\xLtwo$ and $\Phi$ in $\mathcal{L}_2^{0,0}$ or $u_0$
in $\xHone$ and $\Phi$ in $\mathcal{L}_2^{0,1}$, and
\begin{equation}\label{ec1m}
\left\{\begin{array}{rl}i \xdif u^{\epsilon,u_0}&=(\Delta
u^{\epsilon,u_0} + \lambda
|u^{\epsilon,u_0}|^{2\sigma}u^{\epsilon,u_0}-i\alpha
u^{\epsilon,u_0})\xdif t + \sqrt{\epsilon}u^{\epsilon,u_0}\circ\xdif
W,\\
u^{\epsilon,u_0}(0)&=u_0
\end{array}\right.
\end{equation}
with $u_0$ in $\xHone$ and $\Phi$ in $\mathcal{L}_{2,\xR}^{0,s}$
where $s>d/2+1$.\@ When the noise is of multiplicative type, we may
write the equation in terms of a It\^o product,
\begin{equation*}
i \xdif u^{\epsilon,u_0}=(\Delta u^{\epsilon,u_0} + \lambda
|u^{\epsilon,u_0}|^{2\sigma}u^{\epsilon,u_0}-i\alpha
u^{\epsilon,u_0} -(i\epsilon/2)u^{\epsilon,u_0}F_{\Phi})\xdif t +
\sqrt{\epsilon}u^{\epsilon,u_0}\xdif W,
\end{equation*}
where $F_{\Phi}(x)=\sum_{j\in\xN}\left(\Phi e_j (x)\right)^2$ for
$x$ in $\xR^d$ and $\left(e_j\right)_{j\in\xN}$ a complete
orthonormal system of $\xLtwo$.\@ We consider mild solutions; for
example the mild solution of \eqref{ec1} satisfy
\begin{equation*}
\begin{array}{rl}
u^{\epsilon,u_0}(t)=&U(t)u_0-i\lambda\int_0^tU(t-s)(|u^{\epsilon,u_0}(s)|^{2\sigma}
u^{\epsilon,u_0}(s)-i\alpha
u^{\epsilon,u_0}(s))d s\\
& -i\sqrt{\epsilon}\int_0^tU(t-s)d W(s),\quad t>0.
\end{array}
\end{equation*}
The Cauchy problems are globally well posed in $\xLtwo$ and $\xHone$
with the same arguments as in \cite{dBD1}.\\
\indent The main tools in this article are the sample paths LDPs for
the solutions of the three Cauchy problems.\@ They are uniform in
the initial data.\@ Unlike in \cite{DG,EG1,EG2}, we use a
Freidlin-Wentzell type formulation of the upper and lower bounds of
the LDPs.\@ Indeed, it seems that the restriction that initial data
be in compact sets in \cite{EG2} is a real limitation for stochastic
NLS equations.\@ The linear Schr\"odinger group is not compact due
to the lack of smoothing effect and to the fact that we work on the
whole space $\xR^d$.\@ This limitation disappears when we work with
the Freidlin-Wentzell type formulation; we may now obtain bounds for
initial data in balls of $\xLtwo$ or $\xHone$ for $\epsilon$ small
enough.\@ It is well known that in metric spaces and for non uniform
LDPs the two formulations are equivalent.\@ A proof is given in the
Annex and we stress, in
the multiplicative case, on the slight differences with the proof of the result in \cite{EG2}.\\
\indent We denote by $\mathbf{S}(u_0,h)$ the skeleton of equation
\eqref{ec1} or \eqref{ec1m}, {\it i.e.} the mild solution of the
controlled equation
\begin{equation*}
\left\{\begin{array}{l}
i\left(\frac{\xdif u}{\xdif t}+\alpha u\right)=\Delta u +\lambda |u|^{2\sigma}u + \Phi h,\\
u(0)=u_0
\end{array}
\right.
\end{equation*}where $u_0$ belongs to $\xLtwo$ or $\xHone$ in the additive case and the mild
solution of
\begin{equation*}
\left\{\begin{array}{l}
i\left(\frac{\xdif u}{\xdif t}+\alpha u\right)=\Delta u +\lambda |u|^{2\sigma}u + u\Phi h,\\
u(0)=u_0
\end{array}
\right.
\end{equation*}
where $u_0$ belongs to $\xHone$ in the multiplicative case.\\
The rate functions of the LDPs are always defined as
\begin{equation*}
I_T^{u_0}(w)=(1/2)\inf_{h\in\xLtwo\left(0,T;\xLtwo\right):\
\mathbf{S}(u_0,h)=w}\int_0^T\|h(s)\|_{\xLtwo}^2ds.\end{equation*} We
denote for $T$ and $a$ positive by
$K_T^{u_0}(a)=\left(I_T^{u_0}\right)^{-1}\left([0,a]\right)$ the
sets
\begin{equation*}
K_T^{u_0}(a)=\left\{w\in \xCn\left([0,T];\xLtwo\right):\
w=\mathbf{S}(u_0,h),\ (1/2)\int_0^T\|h(s)\|_{\xLtwo}^2ds\leq
a\right\}.
\end{equation*}
We also denote by $d_{\xCn\left([0,T];\xLtwo\right)}$ the usual
distance between sets of $\xCn\left([0,T];\xLtwo\right)$ and by
$d_{\xCn\left([0,T];\xHone\right)}$ the distance between sets of
$\xCn\left([0,T];\xHone\right)$.\\
\indent We write $\tilde{\mathbf{S}}(u_0,f)$ for the skeleton of
equation \eqref{ec1m} where we replace $\Phi h$ by $\frac{\p f}{\p
t}$ where $f$ belongs to $\xHn^1_0\left(0,T;\xHn_{\xR}^{s}\right)$,
the subspace of $\xCn\left([0,T];\xHn_{\xR}^{s}\right)$ of functions
that vanishes at zero and whose time derivative is square
integrable.\@ Also $C_a$ denotes the set
\begin{equation*}
C_a=\left\{f\in\xHn_0^1\left(0,T;\xHn_{\xR}^{s}\right):\ \frac{\p
f}{\p t}\in\xim\Phi,\
I_T^W(f)=(1/2)\left\|\Phi^{-1}_{|(\xker\Phi)^{\perp}}\frac{\p f}{\p
t}\right\|_{\xLtwo(0,T;\xLtwo)}^2\leq a\right\}
\end{equation*}
and $\mathcal{A}(d)$ the set $[2,\infty)$ when $d=1$ or $d=2$ and
$\left[2,2(3d-1)/(3(d-1))\right)$ when $d\geq3$.\@ The above $I_T^W$
is the
good rate function of the LDP for the Wiener process.\\
The uniform LDP with the Freidlin-Wentzell formulation that we need
in the remaining is then as follows.\@ In the additive case we
consider the $\xLtwo$ and $\xHone$ topologies while in the
multiplicative case we consider the $\xHone$ topology only.\@ As it
has been explained previously we do not consider the $\xLtwo$
topology for multiplicative noises since then the $\xLtwo$ norm
remains invariant for the stochastic evolution.\@
\begin{thm}\label{t1}In the additive case and in $\xLtwo$ we have:\\
for every $a$, $\rho$, $T$, $\delta$ and $\gamma$ positive,\\
\begin{tabular}{ll}
(i)&there exists $\epsilon_0$ positive such that for every
$\epsilon$ in $(0,\epsilon_0)$, $u_0$ such that\\
&$\|u_0\|_{\xLtwo}\leq \rho$ and $\tilde{a}$ in $(0,a]$,
\end{tabular}
\begin{equation*}
\mathbb{P}\left(d_{\xCn\left([0,T];\xLtwo\right)}\left(u^{\epsilon,u_0},
K_T^{u_0}(\tilde{a})\right)\geq\delta\right)<\exp\left(-(\tilde{a}
-\gamma)/\epsilon\right),
\end{equation*}
\begin{tabular}{ll}
(ii)&there exists $\epsilon_0$ positive such that for every
$\epsilon$ in $(0,\epsilon_0)$, $u_0$ such that\\
&$\|u_0\|_{\xLtwo}\leq \rho$ and $w$ in $K_T^{u_0}(a)$,
\end{tabular}
\begin{equation*}
\mathbb{P}\left(\left\|u^{\epsilon,u_0}-w
\right\|_{\xCn\left([0,T];\xLtwo\right)}
<\delta\right)>\exp\left(-(I_T^{u_0}(w) +\gamma)/\epsilon\right).
\end{equation*}
In $\xHone$, the result holds for additive and multiplicative noises
replacing in the above $\|u_0\|_{\xLtwo}$ by $\|u_0\|_{\xHone}$ and
$\xCn\left([0,T];\xLtwo\right)$ by
$\xCn\left([0,T];\xHone\right)$.\@
\end{thm}
The proof of this result is given in the annex.\@
\begin{rmrk}The extra condition "For every $a$ positive and $K$ compact in $\xLtwo$, the set
$K_T^{K}(a)=\bigcup_{u_0\in K}K_T^{u_0}(a)$ is a compact subset of
$\xCn\left([0,T];\xLtwo\right)$" often appears to be part of a
uniform LDP.\@ It is not used in the following.\end{rmrk}

\section{Exit from a domain of attraction in $\xLtwo$}\label{s3}
\subsection{Statement of the results}\label{s31}
In this section we only consider the case of an additive noise.\@
Recall that for the real multiplicative noise the mass is decreasing
and thus exit is impossible.\\
\indent We may easily check that the mass
$\mathbf{N}\left(\mathbf{S}(u_0,0)\right)$ of the solution of the
deterministic equation satisfies
\begin{equation}\label{eUA}
\mathbf{N}\left(\mathbf{S}(u_0,0)(t)\right)=
\mathbf{N}\left(u_0\right)\exp\left(-2\alpha t\right).
\end{equation} With noise though, the mass fluctuates around
the deterministic decay.\@ Recall how the It\^o formula applies to
the fluctuation of the mass, see \cite{dBD1} for a proof,
\begin{equation}\label{eI}
\begin{array}{rl}
\mathbf{N}\left(u^{\epsilon,u_0}(t)\right)
-\mathbf{N}\left(u_0\right)=&-2\sqrt{\epsilon}
\pim\int_{\xR^d}\int_0^t\overline{u}^{\epsilon,u_0}dWdx\\
&-2\alpha\left\|u^{\epsilon,u_0}
\right\|_{\xLtwo\left(0,t;\xLtwo\right)}^2 +\epsilon
t\|\Phi\|_{\mathcal{L}_2^{0,0}}^2.\end{array}
\end{equation}
We consider domains $D$ which are bounded measurable subsets of
$\xLtwo$ containing 0 in its interior and invariant by the
deterministic flow, {\it i.e.}
\begin{equation*}
\forall u_0\in D,\ \forall t\geq0,\ \mathbf{S}(u_0,0)(t)\in D.
\end{equation*}
It is thus possible to consider balls.\@ There exists $R$ positive
such
that $D\subset B_R$.\\
We define by
\begin{equation*}
\tau^{\epsilon,u_0}=\inf\left\{t\geq0:\ u^{\epsilon,u_0}(t)\in
D^c\right\}
\end{equation*}
the first exit time of the process $u^{\epsilon,u_0}$ off the domain
$D$.\\
\indent An easy information on the exit time is obtained as
follows.\@ The expectation of an integration via the Duhamel formula
of the It\^o decomposition, the process $u^{\epsilon,u_0}$ being
stopped at the first exit time, gives
$\mathbb{E}\left[\exp\left(-2\alpha\tau^{\epsilon,u_0}\right)\right]=1-
2\alpha R/\left(\epsilon\|\Phi\|_{\mathcal{L}_2^{0,0}}^2\right)$.\@
Without damping we obtain
$\mathbb{E}\left[\tau^{\epsilon,u_0}\right]=
R/\left(\epsilon\|\Phi\|_{\mathcal{L}_2^{0,0}}^2\right)$.\@ To get
more precise
information for small noises we use LDP techniques.\\
\indent Let us introduce
\begin{equation*}
\overline{e}=\inf\left\{I_T^{0}(w):\ w(T)\in \overline{D}^c,\
T>0\right\}.
\end{equation*}
When $\rho$ is positive and small enough, we set
\begin{equation*}
e_{\rho}=\inf\left\{I_T^{u_0}(w):\|u_0\|_{\xLtwo}\leq \rho,\ w(T)\in
\left(D_{-\rho}\right)^c,\ T>0\right\},
\end{equation*}
where $D_{-\rho}=D\setminus\mathcal{N}^0\left(\partial
D,\rho\right)$ and $\partial D$ is the the boundary of $\partial D$
in $\xLtwo$.\@ We define then
\begin{equation*}
\underline{e}=\lim_{\rho\rightarrow0}e_{\rho}.
\end{equation*}
We shall denote in this section by $\|\Phi\|_c$ the norm of $\Phi$
as a bounded operator on $\xLtwo$.\@ Let us start with the following
lemma.\@
\begin{lmm}\label{l0}
$0<\underline{e}\leq\overline{e}$.\@
\end{lmm}
{\bf Proof.}  It is clear that $\underline{e}\leq\overline{e}$.\@
Let us check that $\underline{e}>0$.\@ Let $d$ denote the positive
distance between 0 and $\partial D$.\@ Take $\rho$ small such that
the distance between $B_{\rho}^0$ and $\left(D_{-\rho}\right)^c$ is
larger than $d/2$.\@ Multiplying the evolution equation by
$-i\overline{\mathbf{S}(u_0,h)}$, taking the real part, integrating
over space and using the Duhamel formula we obtain
\begin{equation*}
\begin{array}{l}
\mathbf{N}\left(\mathbf{S}(u_0,h)(T)\right)-\exp\left(-2\alpha
T\right)\mathbf{N}\left(u_0\right)\\=2\int_0^T\exp\left(-2\alpha(T-s)\right)
\pim\left(\int_{\xR^d}\overline{\mathbf{S}\left(u_0,h\right)}\Phi
hdxds\right).
\end{array}
\end{equation*}
If $\mathbf{S}(u_0,h)(T)\in \left(D_{-\rho}\right)^c$ and correspond
to the first escape off $D$ then
\begin{equation*}
\begin{array}{rl}
d/2&\leq2\|\Phi\|_c\int_0^T\exp\left(-2\alpha(T-s)\right)
\left\|\mathbf{S}(u_0,h)(s)\right\|_{\xLtwo}\|h(s)\|_{\xLtwo}ds\\
&\leq2R\|\Phi\|_c\left(\int_0^T\exp\left(-4\alpha(T-s)\right)ds
\right)^{1/2}\|h\|_{\xLtwo(0,T;\xLtwo)},
\end{array}\end{equation*}
thus
\begin{equation*}
\alpha
d^2/\left(8R^2\|\Phi\|_c^2\right)\leq\|h\|_{\xLtwo(0,T;\xLtwo)}^2/2,
\end{equation*}
and the result follows.\hfill $\square$
\begin{rmrk}
We would expect $\underline{e}$ and $\overline{e}$ to be equal.\@ We
may check that it is enough to prove approximate controllability.\@
The argument is however difficult since we are dealing with noises
which are colored space wise, the Schr\"odinger group does not have
global smoothing properties and because of the nonlinearity.\@ If
these two bounds were indeed equal, they would also correspond to
\begin{equation*}
\begin{array}{rl}
\mathcal{E}(D)=&(1/2)\inf\left\{\|h\|_{\xLtwo(0,\infty;\xLtwo)}^2:\
\exists T>0:\ \mathbf{S}(0,h)(T)\in\partial D\right\}\\
=&\inf_{v\in\partial D} V(0,v)
\end{array}
\end{equation*}
where the quasi-potential is defined as
\begin{equation*}
V(u_0,u_f)=\inf\left\{I^{u_0}_T(w):\
w\in\xCn\left([0,\infty);\xLtwo\right),\ w(0)=u_0,\ w(T)=u_f,\
T>0\right\}.
\end{equation*}
\end{rmrk}
\indent We prove in this section the two following results.\@ The
first theorem characterizes the first exit time from the domain.\@
\begin{thm}\label{t2}
For every $u_0$ in $D$ and $\delta$ positive, there exists $L$
positive such that
\begin{equation}\label{et21}
\overline{\lim}_{\epsilon\rightarrow0}\epsilon\log\mathbb{P}\left(
\tau^{\epsilon,u_0}\notin\left(\exp\left((\underline{e}-\delta)/\epsilon\right),
\exp\left((\overline{e}+\delta)/\epsilon\right)\right)\right)\leq-L,
\end{equation}
and for every $u_0$ in $D$,
\begin{equation}\label{et22}
\underline{e}\leq\underline{\lim}_{\epsilon\rightarrow0}\epsilon\log
\mathbb{E}\left(\tau^{\epsilon,u_0}\right)\leq\overline{\lim}_{\epsilon\rightarrow0}
\epsilon\log
\mathbb{E}\left(\tau^{\epsilon,u_0}\right)\leq\overline{e}.
\end{equation}
Moreover, for every $\delta$ positive, there exists $L$ positive
such that
\begin{equation}\label{et23}
\overline{\lim}_{\epsilon\rightarrow0}\epsilon\log\sup_{u_0\in
D}\mathbb{P}\left(\tau^{\epsilon,u_0}\geq\exp\left((\overline{e}
+\delta)/\epsilon\right)\right)\leq-L,
\end{equation}
and
\begin{equation}\label{et24}
\overline{\lim}_{\epsilon\rightarrow0}\epsilon\log\sup_{u_0\in D}
\mathbb{E}\left(\tau^{\epsilon,u_0}\right)\leq\overline{e}.
\end{equation}
\end{thm}
The second theorem characterizes formally the exit points.\@ We
shall define for $\rho$ positive small enough, $N$ a closed subset
of $\partial D$
\begin{equation*}
e_{N,\rho}=\inf\left\{I_T^{u_0}(w):\|u_0\|_{\xLtwo}\leq \rho,\
w(T)\in \left(D\setminus \mathcal{N}^0\left(N,\rho\right)\right)^c,\
T>0\right\}.
\end{equation*}
We then define
\begin{equation*}
\underline{e}_N=\lim_{\rho\rightarrow0}e_{N,\rho}.
\end{equation*}
Note that $e_{\rho}\leq e_{N,\rho}$ and thus
$\underline{e}\leq\underline{e}_N$.\@
\begin{thm}\label{t3} If
$\underline{e}_N>\overline{e}$, then for every $u_0$ in $D$, there
exists $L$ positive such that
\begin{equation*}
\overline{\lim}_{\epsilon\rightarrow
0}\epsilon\log\mathbb{P}\left(u^{\epsilon,u_0}\left(\tau^{\epsilon,u_0}\right)\in
N\right)\leq-L.
\end{equation*}
\end{thm}
Thus the probability of an escape off $D$ via points of $N$ such
that $e_{\rho}\leq e_{N,\rho}$ goes to zero exponentially fast with
$\epsilon$.\\
\indent Suppose that we are able to solve the previous control
problem, then as the noise goes to zero, the probability of an exit
via closed subsets of $\partial D$ where the quasi-potential is not
minimal goes to zero.\@ As the expected exit time is finite, an exit
occurs almost surely.\@ It is exponentially more likely that it
occurs via infima of the quasi-potential.\@ When there are several
infima, the exit measure is a probability measure on $\partial D$.\@
When there is only one infimum we may state the following
corollary.\@
\begin{crllr}
Assume that $v^*$ in $\partial D$ is such that for every $\delta$
positive and $N=\left\{v\in\partial D:\ \|v-v^*\|_{\xLtwo}\geq
\delta\right\}$ we have $\underline{e}_N>\overline{e}$ then
\begin{equation*}
\forall \delta>0,\ \forall u_0\in D,\ \exists L>0:\
\overline{\lim}_{\epsilon\rightarrow0}
\epsilon\log\mathbb{P}\left(\left\|u^{\epsilon,u_0}\left(\tau^{\epsilon,u_0}\right)
-v^*\right\|_{\xLtwo}\geq\delta\right)\leq-L.
\end{equation*}
\end{crllr}

\subsection{Preliminary lemmas}
Let us define
\begin{equation*}
\sigma_{\rho}^{\epsilon,u_0}=\inf\left\{t\geq
0:u^{\epsilon,u_0}(t)\in B_{\rho}^0\cup D^c\right\},
\end{equation*}
where $B_{\rho}^0\subset D$.\@
\begin{lmm}\label{l1}
For every $\rho$ and $L$ positive with $B_{\rho}^0\subset D$, there
exists $T$ and $\epsilon_0$ positive such that for every $u_0$ in
$D$ and $\epsilon$ in $(0,\epsilon_0)$,
\begin{equation*}
\mathbb{P}\left(\sigma_{\rho}^{\epsilon,u_0}>T\right)\leq\exp\left(-L/\epsilon\right).
\end{equation*}
\end{lmm}
{\bf Proof.} The result is straightforward if $u_0$ belongs to
$B_{\rho}^0$.\@ Suppose now that $u_0$ belongs to $D\setminus
B_{\rho}^0$.\@ From equation \eqref{eUA}, the bounded subsets of
$\xLtwo$ are uniformly attracted to zero by the flow of the
deterministic equation.\@ Thus there exists a positive time $T_1$
such that for every $u_1$ in the $\rho/8-$neighborhood of
$D\setminus B_{\rho}^0$ and $t\geq T_1$,
$\mathbf{S}\left(u_1,0\right)(t)\in B_{\rho/8}^0$.\@ We shall choose
$\rho<8$ and follow three steps.\vspace{0.2cm}

\indent{\bf Step 1:} Let us first recall why there exists
$M'=M'(T_1,R,\sigma,\alpha)$ such that
\begin{equation}\label{uni}
\sup_{u_1\in \mathcal{N}^0\left(D\setminus
B_{\rho}^0,\rho/8\right)}\left\|\mathbf{S}(u_1,0)\right\|_{Y^{\left(T_1,2\sigma+2\right)}}\leq
M'.
\end{equation}
From the Strichartz inequalities, there exists $C$ positive such
that
\begin{equation*}
\begin{array}{rl}
\left\|\mathbf{S}(u_1,0)\right\|_{Y^{\left(t,2\sigma+2\right)}}\leq&
C\left\|u_1\right\|_{\xLtwo}+C\left\|\left|\mathbf{S}(u_1,0)\right|^{2\sigma+1}
\right\|_{\xLn^{\gamma'}(0,t;\xLn^{s'})}\\
&+C\alpha\left\|\mathbf{S}\left(u_1,0
\right)\right\|_{\xLn^1\left(0,t;\xLtwo\right)}
\end{array}
\end{equation*}
where $\gamma'$ and $s'$ are such that $1/\gamma'+1/r(\tilde{p})=1$
and $1/s'+1/\tilde{p}=1$ and $\left(r(\tilde{p}),\tilde{p}\right)$
is an admissible pair.\@ Note that the first term is smaller than
$C(R+1)$.\@ From the H\"older inequality, setting
\begin{equation*}
\frac{2\sigma}{2\sigma+2}+\frac{1}{2\sigma+2}=\frac{1}{s'},\ \
\frac{2\sigma}{\omega}+\frac{1}{r(2\sigma+2)}=\frac{1}{\gamma'},
\end{equation*}
we can write
\begin{equation*}
\left\|\left|\mathbf{S}(u_1,0)\right|^{2\sigma+1}\right\|_{\xLn^{\gamma'}
(0,t;\xLn^{s'})}\leq
C\left\|\mathbf{S}(u_1,0)\right\|_{\xLn^{r(2\sigma+2)}(0,t;\xLn^{2\sigma+2})}
\left\|\mathbf{S}(u_1,0)
\right\|_{\xLn^{\omega}(0,t;\xLn^{2\sigma+2})}^{2\sigma}.
\end{equation*}
It is easy to check that since $\sigma<2/d$, we have $\omega<
r(2\sigma+2)$.\@ Thus it follows that
\begin{equation*}
\left\|\mathbf{S}(u_1,0)\right\|_{Y^{\left(t,2\sigma+2\right)}}\leq
C(R+1)+Ct^{\frac{\omega
r(2\sigma+2)}{r(2\sigma+2)-\omega}}\left\|\mathbf{S}(u_1,0)\right\|_{Y^{(t,2\sigma+2)}}^{2\sigma+1}
+C\alpha\sqrt{t}\left\|\mathbf{S}(u_1,0)\right\|_{Y^{(t,2\sigma+2)}}.
\end{equation*}
The function $x\mapsto C(R+1)+Ct^{\frac{\omega
r(2\sigma+2)}{r(2\sigma+2)-\omega}}x^{2\sigma+1}+C\alpha\sqrt{t}x-x$
is positive on a neighborhood of zero.\@ For
$t_0=t_0(R,\sigma,\alpha)$ small enough, the function has at least
one zero.\@ Also, the function goes to $\infty$ as $x$ goes to
$\infty$.\@ Thus, denoting by $M(R,\sigma)$ the first zero of the
above function, we obtain by a classical argument that
$\left\|\mathbf{S}(u_1,0)\right\|_{Y^{\left(t_0,2\sigma+2\right)}}\leq
M(R,\sigma)$
for every $u_1$ in $\mathcal{N}^0\left(D\setminus B_{\rho}^0,\rho/8\right)$.\\
Also, as for every $t$ in $[0,T]$, $\mathbf{S}(u_1,0)(t)$ belongs to
$\mathcal{N}^0\left(D\setminus B_{\rho}^0,\rho/8\right)$, repeating
the previous argument, $u_1$ is replaced by $\mathbf{S}(u_1,0)(t_0)$
and so on, we obtain
\begin{equation*}
\sup_{u_1\in \mathcal{N}^0\left(D\setminus
B_{\rho}^0,\rho/8\right)}\left\|\mathbf{S}(u_1,0)
\right\|_{Y^{\left(T_1,p\right)}}\leq M',
\end{equation*}
where $M'=\left\lceil T_1/t_0\right\rceil M$ proving
\eqref{uni}.\vspace{0.3cm}

\indent{\bf Step 2:} Let us now prove that for $T$ large enough, to
be defined later, and larger than $T_1$, we have
\begin{equation}\label{inclusion}
\mathcal{T}_{\rho}=\left\{w\in\xCn\left([0,T];\xLtwo\right):\
\forall t\in[0,T],\ w(t)\in\mathcal{N}^0\left(D\setminus
B_{\rho}^0,\rho/8\right)\right\}\subset K_T^{u_0}(2L)^c.
\end{equation}
Since $K_T^{u_0}(2L)$ is included in the image of
$\mathbf{S}(u_0,\cdot)$ it suffices to consider $w$ in
$\mathcal{T}_{\rho}$ such that $w=\mathbf{S}(u_0,h)$ for some $h$ in
$\xLtwo(0,T;\xLtwo)$.\@ Take $h$ such that $\mathbf{S}(u_0,h)$
belongs to $\mathcal{T}_{\rho}$ we have
\begin{equation*}
\left\|\mathbf{S}(u_0,h)-\mathbf{S}(u_0,0)\right\|_{\xCn\left([0,T_1];\xLtwo\right)}
\geq\left\|\mathbf{S}(u_0,h)(T_1)-\mathbf{S}(u_0,0)(T_1)\right\|_{\xLtwo}\geq3\rho/4,
\end{equation*}
but also, necessarily, for the admissible pair
$\left(r(2\sigma+2),2\sigma+2\right)$,
\begin{equation}\label{e5}
\left\|\mathbf{S}(u_0,h)-\mathbf{S}(u_0,0)\right\|_{Y^{\left(T_1,2\sigma+2\right)}}
\geq3\rho/4.
\end{equation}
Denote by $\mathbf{S}^{M'+1}$ the skeleton corresponding to the
following control problem
\begin{equation*}
\left\{\begin{array}{l} i\left(\frac{\xdif u}{\xdif t}+\alpha
u\right)=\Delta u +\lambda \theta\left(
\frac{\left\|u\right\|_{Y^{(t,2\sigma+2)}}}{M'+1}\right)|u|^{2\sigma}u + \Phi h,\\
u(0)=u_1
\end{array}
\right.
\end{equation*}
where $\theta$ is a $\xCn^{\infty}$ function with compact support,
such that $\theta(x)=0$ if $x\geq2$ and $\theta(x)=1$ if $0\leq
x\leq1$.\@ Then \eqref{e5} implies that
\begin{equation*}
\left\|\mathbf{S}^{M'+1}(u_0,h)-\mathbf{S}^{M'+1}(u_0,0)\right\|_{Y^{\left(T_1,2\sigma+2\right)}}
\geq3\rho/4.
\end{equation*}
We shall now split the interval $[0,T_1]$ in many parts.\@ We shall
denote here by $Y^{s,t,2\sigma+2}$ for $s<t$ the space
$Y^{t,2\sigma+2}$ on the interval $[s,t]$.\@ Applying the Strichartz
inequalities on a small interval $[0,t]$ with the computations in
the proof of Lemma 3.3 in \cite{dBD0}, we obtain
\begin{equation*}
\begin{array}{r}
\left\|\mathbf{S}^{M'+1}(u_0,h)-\mathbf{S}^{M'+1}(u_0,0)
\right\|_{Y^{\left(t,2\sigma+2\right)}}\leq
C\alpha\sqrt{t}\left\|\mathbf{S}^{M'+1}(u_0,h)-\mathbf{S}^{M'+1}
(u_0,0)\right\|_{Y^{\left(t,2\sigma+2\right)}}\\
+C_{M'+1}t^{1-d\sigma/2}\left\|\mathbf{S}^{M'+1}(u_0,h)
-\mathbf{S}^{M'+1}(u_0,0) \right\|_{Y^{\left(t,2\sigma+2\right)}}
+C\sqrt{t}\|\Phi\|_c\|h\|_{\xLtwo(0,t;\xLtwo)}
\end{array}
\end{equation*}
where $C_{M'+1}$ is a constant which depends on $M'+1$.\@ Take $t_1$
small enough such that $C_{M'+1}t_1^{1-d\sigma
d/2}+C\alpha\sqrt{t_1}\leq1/2$.\@ We obtain then
\begin{equation*}
\left\|\mathbf{S}^{M'+1}(u_0,h)-\mathbf{S}^{M'+1}(u_0,0)
\right\|_{Y^{\left(t_1,2\sigma+2\right)}}\leq
2C\sqrt{t_1}\|\Phi\|_c\|h\|_{\xLtwo(0,t_1;\xLtwo)}.
\end{equation*}
In the case where $2t_1<T_1$, let us see how such inequality
propagates on $\left[t_1,2t_1\right]$.\@ We now have two different
initial data $\mathbf{S}^{M'+1}(u_0,h)\left(t_1\right)$ and
$\mathbf{S}^{M'+1}(u_0,0)\left(t_1\right)$.\@ We obtain similarly
\begin{equation*}
\begin{array}{l}
\left\|\mathbf{S}^{M'+1}(u_0,h)-\mathbf{S}^{M'+1}(u_0,0)
\right\|_{Y^{\left(t_1,2t_1,2\sigma+2\right)}}\\
\leq
2C\sqrt{t_1}\|\Phi\|_c\|h\|_{\xLtwo(0,t_1;\xLtwo)}+2\left\|\mathbf{S}^{M'+1}(u_0,h)\left(t_1\right)
-\mathbf{S}^{M'+1}(u_0,0)\left(t_1\right)\right\|_{\xHone}\\
\leq
2C\sqrt{t_1}\|\Phi\|_c\|h\|_{\xLtwo(0,T_1;\xLtwo)}+2\left\|\mathbf{S}^{M'+1}(u_0,h)
\left(t_1\right)-\mathbf{S}^{M'+1}(u_0,0)\left(t_1\right)\right\|_{Y^{\left(0,t_1,2\sigma+2\right)}}.
\end{array}
\end{equation*}
Then iterating on each interval of the form $[kt_1,(k+1)t_1]$ for
$k$ in
$\left\{1,...,\left\lfloor\left.T_1/t_1-1\right\rfloor\right.\right\}$,
the remaining term can be treated similarly, and using the triangle
inequality we obtain that
\begin{equation*}
\left\|\mathbf{S}^{M'+1}(u_0,h)-\mathbf{S}^{M'+1}(u_0,0)
\right\|_{Y^{\left(T_1,2\sigma+2\right)}}\leq
2^{\left\lceil\left.T_1/t_1\right\rceil\right.+1}
C\sqrt{t_1}\|\Phi\|_c \|h\|_{\xLtwo(0,t_1;\xLtwo)}.
\end{equation*}
We may then conclude that
\begin{equation*}
\left\|h\right\|_{\xLtwo\left(0,T_1;\xLtwo\right)}^2/2\geq M''
\end{equation*}
where $M''=\rho^2/(8C\left(t_1,T_1\right)\|\Phi\|_c^2)$ and
$C\left(t_1,T_1\right)$ is a constant which depends only on $t_1$
and $T_1$.\@ Note that we have used for later purposes that
$3\rho/2>\rho/2$.\vspace{0.3cm}

\indent Similarly replacing $[0,T_1]$ by $[T_1,2T_1]$ and $u_0$
respectively by $\mathbf{S}\left(u_0,h\right)(T_1)$ and
$\mathbf{S}\left(u_0,0\right)(T_1)$ in \eqref{e5}, the inequality
still holds true.\@ Thus thanks to the inverse triangle inequality
we obtain on $[T_1,2T_1]$
\begin{equation*}
\begin{array}{l}
\left\|\mathbf{S}^{M'+1}(u_0,h)-\mathbf{S}^{M'+1}(u_0,0)\right\|_{Y^{\left(T_1,2T_1,2\sigma+2\right)}}\\
=\left\|\mathbf{S}^{M'+1}\left(\mathbf{S}^{M'+1}(u_0,h)(T_1),h\right)-\mathbf{S}^{M'+1}\left(\mathbf{S}^{M'+1}(u_0,0)(T_1),0\right)\right\|_{Y^{\left(0,T_1,2\sigma+2\right)}}\\
\geq3\rho/4
\end{array}
\end{equation*}
Thus from the inverse triangle inequality along with the fact that
for both $\mathbf{S}^{M'+1}(u_0,h)(T_1)$ and
$\mathbf{S}^{M'+1}(u_0,0)(T_1)$ as initial data the deterministic
solutions belong to the ball $B_{\rho/8}^0$, we obtain
\begin{equation*}
\left\|\mathbf{S}^{M'+1}\left(\mathbf{S}^{M'+1}(u_0,h)(T_1),h\right)-\mathbf{S}^{M'+1}\left(\mathbf{S}^{M'+1}(u_0,h)(T_1),0\right)\right\|_{Y^{\left(0,T_1,2\sigma+2\right)}}\\
\geq\rho/2.
\end{equation*}
We finally obtain the same lower bound
\begin{equation*}
\left\|h\right\|_{\xLtwo\left(T_1,2T_1;\xLtwo\right)}^2/2\geq M''
\end{equation*}
as before.\\
Iterating the argument we obtain if $T>2T_1$,
\begin{equation*}
\left\|h\right\|_{\xLtwo\left(0,2T_1;\xLtwo\right)}^2/2
=\left\|h\right\|_{\xLtwo\left(0,T_1;\xLtwo\right)}^2/2
+\left\|h\right\|_{\xLtwo\left(T_1,2T_1;\xLtwo\right)}^2/2\geq 2M''.
\end{equation*}
Thus for $j$ positive and $T>jT_1$, we obtain, iterating the above
argument, that
\begin{equation*}
\left\|h\right\|_{\xLtwo\left(0,jT_1;\xLtwo\right)}^2/2\geq jM''.
\end{equation*}
The result \eqref{inclusion} is obtained for $T=jT_1$ where $j$ is
such that $jM''>2L$.\vspace{0.3cm}

\indent{\bf Step 3:} We may now conclude from the (i) of Theorem
\ref{t1} since,
\begin{equation*}
\begin{array}{rl}
\mathbb{P}\left(\sigma_{\rho}^{\epsilon,u_0}>T\right)&=\mathbb{P}
\left(\forall t\in[0,T],\ u^{\epsilon,u_0}(t)\in D\setminus B_{\rho}^0\right)\\
&=\mathbb{P}\left(d_{\xCn\left([0,T];\xLtwo\right)}
\left(u^{\epsilon,u_0},\mathcal{T}_{\rho}^c\right)>\rho/8\right),\\
&\leq\mathbb{P}\left(d_{\xCn\left([0,T];\xLtwo\right)}
\left(u^{\epsilon,u_0},K_T^{u_0}(2L)\right)\geq\rho/8\right),
\end{array}
\end{equation*}
taking $a=2L$, $\rho=R$ where $D\subset B_R$,
$\delta=\rho/8$ and $\gamma=L$.\\
\indent Note that if $\rho\geq 8$, we should replace $R+1$ by
$R+\rho/8$ and $M'+1$ by $M'+\rho/8$.\@ Anyway, we will use the
lemma for small $\rho$.\hfill$\square$
\begin{lmm}\label{l2}
For every $\rho$ positive such that $B_{\rho}^0\subset D$ and $u_0$
in $D$, there exists $L$ positive such that
\begin{equation*}
\overline{\lim}_{\epsilon\rightarrow 0}\epsilon\log
\mathbb{P}\left(u^{\epsilon,u_0}\left(\sigma_{\rho}^{\epsilon,u_0}
\right)\in\partial D \right)\leq-L
\end{equation*}
\end{lmm}
{\bf Proof.}  Take $\rho$ positive satisfying the assumptions of the
lemma and take $u_0$ in $D$.\@ When $u_0$ belongs to $B_{\rho}^0$
the result is straightforward.\@ Suppose now that $u_0$ belongs to
$D\setminus B_{\rho}^0$.\@ Let $T$ be defined as
\begin{equation*}
T=\inf\left\{t\geq 0:\ \mathbf{S}\left(u_0,0\right)(t)\in
B_{\rho/2}^0\right\},
\end{equation*}
then since $\mathbf{S}\left(u_0,0\right)\left([0,T]\right)$ is a
compact subset of $D$, the distance $d$ between
$\mathbf{S}\left(u_0,0\right)\left([0,T]\right)$ and $D^c$ is well
defined and positive.\@ The conclusion follows then from the fact
that
\begin{equation*}
\mathbb{P}\left(u^{\epsilon,u_0}\left(\sigma_{\rho}^{\epsilon,u_0}\right)\in\partial
D\right)\leq \mathbb{P}\left(\left\|u^{\epsilon,u_0}-\mathbf{S}
\left(u_0,0\right)\right\|_{\xCn\left([0,T];\xLtwo\right)}\geq(\rho\wedge
d)/2\right),
\end{equation*}
the LDP and the fact that, from the compactness of the sets
$K_T^{u_0}(a)$ for $a$ positive, we have
\begin{equation*}
\inf_{h\in\xLtwo\left(0,T;\xLtwo\right):\
\|\mathbf{S}(u_0,h)-\mathbf{S}(u_0,0)\|_{\xCn\left([0,T];\xLtwo\right)}\geq(\rho\wedge
d)/2}\|h\|_{\xLtwo(0,T;\xLtwo)}^2>0.
\end{equation*}
We have used the fact that the upper bound of the LDP in the
Freidlin-Wentzell formulation implies the classical upper bound.\@
Note that this is a well known result for non uniform LDPs.\@ Indeed
we do not need a uniform LDP in this proof. \hfill
$\square$\vspace{0.3cm}

\indent The following lemma replaces Lemma 5.7.23 in \cite{DZ}.\@
Indeed, the case of a stochastic PDE is more intricate than that of
a SDE since the linear group is only strongly and not uniformly
continuous.\@ However, it is possible to prove that the group on
$\xLtwo$ when acting on bounded sets of $\xHone$ is uniformly
continuous.\@ We shall proceed in a different manner and thus we do
not loose in regularity.\@
\begin{lmm}\label{l3} For every
$\rho$ and $L$ positive such that $B_{2\rho}^0\subset D$, there
exists $T(L,\rho)<\infty$ such that
\begin{equation*}
\overline{\lim}_{\epsilon\rightarrow0}\epsilon\log\sup_{u_0\in
S_{\rho}^0}\mathbb{P}\left(\sup_{t\in[0,T(L,\rho)]}\left(\mathbf{N}
\left(u^{\epsilon,u_0}(t)\right)-\mathbf{N}\left(u_0\right)\right)\geq
3\rho^2\right)\leq-L
\end{equation*}
\end{lmm}
{\bf Proof.} Take $L$ and $\rho$ positive.\@ Note that for every
$\epsilon$ in $(0,\epsilon_0)$ where
$\epsilon_0=\rho^2/\|\Phi\|_{\mathcal{L}_2^{0,0}}^2$, for
$T(L,\rho)\leq 1$ we have $\epsilon
T(L,\rho)\|\Phi\|_{\mathcal{L}_2^{0,0}}^2<\rho^2$.\@ Thus from
equation \eqref{eI}, we know that it is enough to prove that there
exists $T(L,\rho)\leq1$ such that for $\epsilon_1$ small enough,
$\epsilon_1<\epsilon_0$, and all $\epsilon<\epsilon_0$,
\begin{equation*}
\epsilon\log\sup_{u_0\in
S_{\rho}^0}\mathbb{P}\left(\sup_{t\in[0,T(L,\rho)]}
\left(-2\sqrt{\epsilon}\pim\int_{\xR^d}\int_0^t\overline{u}^{\epsilon,u_0,\tau}dWdx\right)\geq
2\rho^2\right)\leq-L,\end{equation*} where $u^{\epsilon,u_0,\tau}$
is the process $u^{\epsilon,u_0}$ stopped at
$\tau_{S_{2\rho}^0}^{\epsilon,u_0}$, the first time when
$u^{\epsilon,u_0}$ hits $S_{2\rho}^0$.\@ Setting
$Z(t)=\pim\int_{\xR^d}\int_0^t\overline{u}^{\epsilon,u_0,\tau}dWdx$,
it is enough to show that
\begin{equation*}
\epsilon\log\sup_{u_0\in
S_{\rho}^0}\mathbb{P}\left(\sup_{t\in[0,T(L,\rho)]}\left|Z(t)\right|\geq
\rho^2/\sqrt{\epsilon}\right)\leq-L,\end{equation*} and thus to show
exponential tail estimates for the process $Z(t)$.\@ Our proof now
follows closely that of \cite{P2}[Theorem 2.1].\@ We introduce the
function $f_l(x)=\sqrt{1+lx^2}$, where $l$ is a positive
parameter.\@ We now apply the It\^o formula to $f_l(Z(t))$ and the
process decomposes into $1+E_l(t)+R_l(t)$ where
\begin{equation*}
E_l(t)=\int_0^t\frac{2lZ(t)}{\sqrt{1+lZ(t)^2}}dZ(t)
-(1/2)\int_0^t\left(\frac{2lZ(t)}{\sqrt{1+lZ(t)^2}}\right)^2d<Z>_t,
\end{equation*}
and
\begin{equation*}
R_l(t)=(1/2)\int_0^t\left(\frac{2lZ(t)}{\sqrt{1+lZ(t)^2}}\right)^2d<Z>_t
+\int_0^t\frac{l}{\left(1+lZ(t)^2\right)^{3/2)}}d<Z>_t.
\end{equation*}
Moreover, given $\left(e_j\right)_{j\in\xN}$ a complete orthonormal
system of $\xLtwo$,
\begin{equation*}<Z(t)>=\int_0^t\sum_{j\in\xN}\left(u^{\epsilon,u_0,\tau},-i\Phi
e_j\right)_{\xLtwo}^2(s)ds, \end{equation*} we prove with the
H\"older inequality that $|R_l(t)|\leq
12l\rho^2\|\Phi\|_{\mathcal{L}^{0,0}_2}^2t$, for every $u_0$ in
$D$.\@ We may thus write
\begin{equation*}
\begin{array}{rl}
&\mathbb{P}\left(\sup_{t\in[0,T(L,\rho)]}\left|Z(t)\right|\geq
\rho^2/\sqrt{\epsilon}\right)\\&=\mathbb{P}
\left(\sup_{t\in[0,T(L,\rho)]}\exp\left(f_l(Z(t))
\right)\geq\exp\left(f_l\left(\rho^2/\sqrt{\epsilon}\right)\right)\right)\\
&\leq\mathbb{P}\left(\sup_{t\in[0,T(L,\rho)]}\exp\left(E_l(t)
\right)\geq\exp\left(f_l\left(\rho^2/\sqrt{\epsilon}\right)
-1-12l\rho^2\|\Phi\|_{\mathcal{L}^{0,0}_2}^2T(L,\rho)\right)\right).
\end{array}
\end{equation*}
The Novikov condition is also satisfied and $E_l(t)$ is such that
$\left(\exp\left(E_l(t)\right)\right)_{t\in\xR^+}$ is a uniformly
integrable martingale.\@ The exponential tail estimates follow from
the Doob inequality optimizing on the parameter $l$.\@ We may then
write
\begin{equation*}
\sup_{u_0\in
S_{\rho}^0}\mathbb{P}\left(\sup_{t\in[0,T(L,\rho)]}\left|Z(t)\right|\geq
\rho^2/\sqrt{\epsilon}\right)\leq3\exp\left(-\frac{\rho^2}{48\epsilon
\|\Phi\|_{\mathcal{L}_2^{0,0}}^2T(L,\rho)}\right).\end{equation*} We
now conclude setting
$T(L,\rho)=\rho^2/\left(50\|\Phi\|_{\mathcal{L}_2^{0,0}}^2L\right)$
and choosing $\epsilon_1<\epsilon_0$ small enough.\hfill $\square$

\subsection{Proof of Theorem \ref{t2} and Theorem \ref{t3}}
We first prove Theorem \ref{t2}.\\
{\bf Proof of Theorem \ref{t2}.} Let us first prove \eqref{et24} and
deduce \eqref{et23}.\@ Fix $\delta$ positive and choose $h$ and
$T_1$ such that $\mathbf{S}(0,h)(T_1)\in\overline{D}^c$ and
\begin{equation*}
I_{T_1}^{0}\left(\mathbf{S}(0,h)\right)
=(1/2)\|h\|_{\xLtwo\left(0,T;\xLtwo\right)}^2\leq
\overline{e}+\delta/5.
\end{equation*}
Let $d_0$ denote the positive distance between
$\mathbf{S}(0,h)\left(T_1\right)$ and $\overline{D}$.\@ With similar
arguments as in \cite{dBD1} or with a truncation argument we may
prove that the skeleton is continuous with respect to the initial
datum for the $\xLtwo$ topology.\@ Thus there exists $\rho$
positive, a function of $h$ which has been fixed, such that if $u_0$
belongs to $B_{\rho}^0$ then
\begin{equation*}
\left\|\mathbf{S}\left(u_0,h\right)-\mathbf{S}(0,h)\right\|_{\xCn\left([0,T_1];\xLtwo\right)}
<d_0/2.
\end{equation*}
We may assume that $\rho$ is such that $B_{\rho}^0\subset D$.\@ From
the triangle inequality and the (ii) of Theorem \ref{t1}, there
exists $\epsilon_1$ positive such that for all $\epsilon$ in
$\left(0,\epsilon_1\right)$ and $u_0$ in $B_{\rho}^0$,
\begin{equation*}
\begin{array}{rl}
\mathbb{P}\left(\tau^{\epsilon,u_0}<T_1\right)&\geq
\mathbb{P}\left(\left\|u^{\epsilon,u_0}-\mathbf{S}(0,h)
\right\|_{\xCn\left([0,T_1];\xLtwo\right)}<d_0\right)\\
&\geq\mathbb{P}\left(\left\|u^{\epsilon,u_0}-\mathbf{S}
\left(u_0,h\right)\right\|_{\xCn\left([0,T_1];\xLtwo\right)}<d_0/2\right)\\
&\geq\exp\left(-(I_{T_1}^{u_0}\left(\mathbf{S}(u_0,h)\right)
+\frac{\delta}{5})/\epsilon\right).
\end{array}
\end{equation*}
From Lemma \ref{l1}, there exists $T_2$ and $\epsilon_2$ positive
such that for all $\epsilon$ in $\left(0,\epsilon_2\right)$,
\begin{equation*}
\inf_{u_0\in D} \mathbb{P}\left(\sigma_{\rho}^{\epsilon,u_0}\leq
T_2\right)\geq 1/2.
\end{equation*}
Thus, for $T=T_1+T_2$, from the strong Markov property we obtain
that for all $\epsilon<\epsilon_3<\epsilon_1\wedge\epsilon_2$.\@
\begin{equation*}
\begin{array}{rl}
q=\inf_{u_0\in D}\mathbb{P}\left(\tau^{\epsilon,u_0}\leq
T\right)&\geq \inf_{u_0\in
D}\mathbb{P}\left(\sigma_{\rho}^{\epsilon,u_0}\leq
T_2\right)\inf_{u_0\in
B_{\rho}^0}\mathbb{P}\left(\tau^{\epsilon,u_0}\leq T_1\right)\\
&\geq(1/2)\exp\left(-\left(I_{T_1}^{u_0}
\left(\mathbf{S}(u_0,h)\right)+\delta/5\right)/\epsilon\right)\\
&\geq \exp\left(-\left(I_{T_1}^{u_0}\left(\mathbf{S}(u_0,h)\right)
+2\delta/5\right)/\epsilon\right).
\end{array}
\end{equation*}
Thus, for any $k\geq1$, we have
\begin{equation*}
\begin{array}{rl}
\mathbb{P}\left(\tau^{\epsilon,u_0}>(k+1)T\right)&=
\left[1-\mathbb{P}\left(\tau^{\epsilon,u_0}\leq(k+1)T|\tau^{\epsilon,u_0}>kT\right)
\right]\mathbb{P}\left(\tau^{\epsilon,u_0}>kT\right)\\
&\leq(1-q)\mathbb{P}\left(\tau^{\epsilon,u_0}>kT\right)\\
&\leq(1-q)^k.
 \end{array}
\end{equation*}
We may now compute, since
$I_{T_1}^{u_0}\left(\mathbf{S}\left(u_0,h\right)\right)=
I_{T_1}^{0}\left(\mathbf{S}\left(0,h\right)\right)=(1/2)\|h\|_{\xLtwo\left(0,T;\xLtwo\right)}^2$
\begin{equation*}
\begin{array}{rl}
\sup_{u_0\in
D}\mathbb{E}\left(\tau^{\epsilon,u_0}\right)&=\sup_{u_0\in
D}\int_0^{\infty}\mathbb{P}\left(\tau^{\epsilon,u_0}>t\right)dt\\
&\leq T\left[1+\sum_{k=1}^{\infty}\sup_{x\in
D}\mathbb{P}\left(\tau^{\epsilon,u_0}>kT\right)\right]\\
&\leq T/q\\
&\leq T\exp\left((\overline{e}+3\delta/5)/\epsilon\right).
\end{array}
\end{equation*}
It implies that there exists $\epsilon_4$ small enough such that for
$\epsilon$ in $\left(0,\epsilon_4\right)$,
\begin{equation}\label{e2}
\sup_{u_0\in D}\mathbb{E}\left(\tau^{\epsilon,u_0}\right)\leq
\exp\left((\overline{e}+4\delta/5)/\epsilon\right).
\end{equation}
Thus the Chebychev inequality gives that
\begin{equation*}
\sup_{u_0\in D}\mathbb{P}\left(\tau^{\epsilon,u_0}\geq\exp
\left((\overline{e}+\delta)/\epsilon\right)\right)
\leq\exp\left(-(\overline{e}+\delta)/\epsilon\right)\sup_{u_0\in
D}\mathbb{E}\left(\tau^{\epsilon,u_0}\right),
\end{equation*}
in other words
\begin{equation}\label{e3} \sup_{u_0\in
D}\mathbb{P}\left(\tau^{\epsilon,u_0}\geq\exp
\left((\overline{e}+\delta)/\epsilon\right)\right)
\leq\exp\left(-\delta/(5\epsilon)\right).
\end{equation}
Relations \eqref{e2} and \eqref{e3} imply \eqref{et24} and
\eqref{et23}.\vspace{0.3cm}

\indent Let us now prove the lower bound on $\tau^{\epsilon,u_0}$.\@
Take $\delta$ positive.\@ Remind that we have proved that
$\underline{e}>0$.\@ Take $\rho$ positive small enough such that
$\underline{e}-\delta/4\leq e_{\rho}$ and $B_{2\rho}^0\subset D$.\@
We define the following sequences of stopping times, $\theta_0=0$
and for $k$ in $\xN$,
\begin{equation*}\begin{array}{rl}
\tau_k&=\inf\left\{t\geq\theta_k:\ u^{\epsilon,u_0}(t)\in
B_{\rho}^0\cup
D^c\right\},\\
\theta_{k+1}&=\inf\left\{t>\tau_k:\ u^{\epsilon,u_0}(t)\in
S_{2\rho}^0\right\},\end{array}
\end{equation*}
where $\theta_{k+1}=\infty$ if $u^{\epsilon,u_0}(\tau_k)\in
\partial D$.\@ Fix $T_1=T\left(\underline{e}-3\delta/4,\rho\right)$
given in Lemma \ref{l3}.\@ We know that there exists $\epsilon_1$
positive such that for all $\epsilon$ in $(0,\epsilon_1)$, for all
$k\geq 1$ and $u_0$ in $D$,
\begin{equation*}
\mathbb{P}\left(\theta_k-\tau_{k-1}\leq
T_1\right)\leq\exp\left(-(\underline{e}-3\delta/4)/\epsilon\right).
\end{equation*}
For $u_0$ in $D$ and an $m$ in $\xN^*$, we have
\begin{equation}\label{e4}
\begin{array}{rl}
\mathbb{P}\left(\tau^{\epsilon,u_0}\leq mT_1\right)\leq&
\mathbb{P}\left(\tau^{\epsilon,u_0}=\tau_0\right)
+\sum_{k=1}^m\mathbb{P}\left(\tau^{\epsilon,u_0}=\tau_k\right)\\
&+\mathbb{P}\left(\exists k\in\{1,...,m\}:\ \theta_k-\tau_{k-1}\leq
T_1\right)\\
=&\mathbb{P}\left(\tau^{\epsilon,u_0}=\tau_0\right)
+\sum_{k=1}^m\mathbb{P}\left(\tau^{\epsilon,u_0}=\tau_k\right)\\
&+\sum_{k=1}^m\mathbb{P}\left(\theta_k-\tau_{k-1}\leq T_1\right).
\end{array}
\end{equation}
In other words the escape before $mT_1$ can occur either as an
escape without passing in the small ball $B_{\rho}^0$ (if $u_0$
belongs to $D\setminus B_{\rho}^0$) or as an escape with $k$ in
$\{1,...m\}$ significant fluctuations off $B_{\rho}^0$, {\it i.e.}
crossing $S_{2\rho}^0$, or at least one of the $m$ first transitions
between $S_{\rho}^0$ and $S_{2\rho}^0$ happens in less than $T_1$.\@
The latter is known to be arbitrarily small.\@
Let us prove that the remaining probabilities are small enough for small $\epsilon$.\\
For every $k\geq 1$ and $T_2$ positive, we may write
\begin{equation*}
\mathbb{P}\left(\tau^{\epsilon,u_0}=\tau_k\right)\leq\mathbb{P}\left(\tau^{\epsilon,u_0}\leq
T_2;\tau^{\epsilon,u_0}=\tau_k\right)+\mathbb{P}\left(\sigma_{\rho}^{\epsilon,u_0}>T_2\right).
\end{equation*}
Fix $T_2$ as in Lemma \ref{l1} with $L=\underline{e}-3\delta/4$.\@
Thus there exists $\epsilon_2$ small enough such that for $\epsilon$
in $(0,\epsilon_2)$,
\begin{equation*}
\mathbb{P}\left(\sigma_{\rho}^{\epsilon,u_0}>T_2\right)\leq
\exp\left(-(\underline{e}-3\delta/4)/\epsilon\right).
\end{equation*}
Also, from the (i) of Theorem \ref{t1}, we obtain that there exists
$\epsilon_3$ positive such that for every $u_1$ in $B_{\rho}^{0}$
and $\epsilon$ in $\left(0,\epsilon_3\right)$,
\begin{equation*}
\begin{array}{rl}
\mathbb{P}\left(\tau^{\epsilon,u_1}\leq
T_2\right)&\leq\mathbb{P}\left(d_{\xCn\left(\left[0,T_2\right];\xLtwo\right)}
\left(u^{\epsilon,u_1},K_{T_2}^{u_1}\left(e_{\rho}
-\delta/4\right)\right)\geq\rho\right)\\
&\leq\exp\left(-(e_{\rho}-\delta/2)/\epsilon\right)\\
&\leq\exp\left(-(\underline{e}-3\delta/4)/\epsilon\right).
\end{array}
\end{equation*}
Thus the above bound holds for
$\mathbb{P}\left(\tau^{\epsilon,u_0}\leq
T_2;\tau^{\epsilon,u_0}=\tau_k\right)$ replacing $u_1$ by
$u^{\epsilon,u_0}\left(\tau_{k-1}\right)$ since as $k\geq1$,
$u^{\epsilon,u_0}\left(\tau_{k-1}\right)$ belongs to $B_{\rho}^0$
and $\tau_k-\tau_{k-1}\leq T_2$ and using the Markov property.\@ The
inequality \eqref{e4} gives that for all $\epsilon$ in
$\left(0,\epsilon_0\right)$ where
$\epsilon_0=\epsilon_1\wedge\epsilon_2\wedge\epsilon_3$,
\begin{equation*}
\mathbb{P}\left(\tau^{\epsilon,u_0}\leq mT_1\right)\leq
\mathbb{P}\left(u^{\epsilon,u_0}\left(\sigma_{\rho}^{\epsilon,u_0}\right)\in\partial
D\right)+3m\exp\left(-(\underline{e}-3\delta/4)/\epsilon\right).
\end{equation*}
Fix $m=\left\lceil
(1/T_1)\exp\left((\underline{e}-\delta)/\epsilon\right)\right\rceil$,
then for all $\epsilon$ in $\left(0,\epsilon_0\right)$,
\begin{equation*}\begin{array}{rl}
\mathbb{P}\left(\tau^{\epsilon,u_0}\leq
\exp\left((\underline{e}-\delta)/\epsilon\right)\right)&\leq
\mathbb{P}\left(\tau^{\epsilon,u_0}\leq mT_1\right)\\
&\leq\mathbb{P}\left(u^{\epsilon,u_0}\left(\sigma_{\rho}^{\epsilon,u_0}\right)\in\partial
D\right)+(3/T_1)\exp\left(-\delta/(4\epsilon)\right).
\end{array}
\end{equation*}
We may now conclude with Lemma \ref{l2} and obtain the expected
lower bound on $\mathbb{E}\left(\tau^{\epsilon,u_0}\right)$ from the
Chebychev inequality.\hfill $\square$\vspace{0.3cm}

Let us now prove Theorem \ref{t3}.\\
{\bf Proof of Theorem \ref{t3}.} Let $N$ be closed subset of
$\partial D$.\@ When $\underline{e}_N=\infty$ we shall replace in
the proof that follows $\underline{e}_N$ by an increasing sequence
of positive numbers.\@ Take $\delta$ such that
$0<\delta<(\underline{e}_N-\overline{e})/3$, $\rho$ positive such
that $\underline{e}_N-\delta/3\leq e_{N,\rho}$ and
$B_{2\rho}^0\subset D$.\@ Define the same sequences of stopping
times $\left(\tau_k\right)_{k\in\xN}$ and
$\left(\theta_k\right)_{k\in\xN}$
as in the proof of Theorem \ref{t2}.\\
Take $L=\underline{e}_N-\delta$ and $T_1$ and
$T_2=T\left(L,\rho\right)$ as in Lemma \ref{l1} and \ref{l3}.\@
Thanks to Lemma \ref{l1} and the uniform LDP, with a computation
similar to the one following inequality \eqref{e4}, we obtain that
for $\epsilon_0$ small enough and $\epsilon\leq \epsilon_0$,
\begin{equation*}
\begin{array}{l}
\sup_{u_0\in S_{2\rho}^0}\mathbb{P}\left(u^{\epsilon,u_0}
\left(\sigma_{\rho}^{\epsilon,u_0}\right)\in N \right)\\
\leq\sup_{u_0\in S_{2\rho}^0}\mathbb{P}\left(u^{\epsilon,u_0}
\left(\sigma_{\rho}^{\epsilon,u_0}\right)\in
N,\sigma_{\rho}^{\epsilon,u_0}\leq T_1\right)
+\sup_{u_0\in S_{2\rho}^0}\mathbb{P}\left(\sigma_{\rho}^{\epsilon,u_0}>T_1\right)\\
\leq\sup_{u_0\in
B_{2\rho}^0}\mathbb{P}\left(d_{\xCn\left([0,T_1];\xLtwo\right)}
\left(u^{\epsilon,u_0},K_{T_1}^{u_0}\left(e_{N,\rho}-\delta/3\right)\right)
\geq\rho\right)
\\
\ \ +\sup_{u_0\in D}\mathbb{P}\left(\sigma_{\rho}^{\epsilon,u_0}>T_1\right)\\
\leq2\exp\left(-(\underline{e}_N-\delta)/\epsilon\right).
\end{array}
\end{equation*}
Possibly choosing $\epsilon_0$ smaller, we may assume that for every
positive integer $l$ and every $\epsilon\leq\epsilon_0$,
\begin{equation*}
\begin{array}{rl}
\sup_{u_0\in D}\mathbb{P}\left(\tau_l\leq lT_2\right)\leq&
l\sup_{u_0\in S_{\rho}^0}
\mathbb{P}\left(\sup_{t\in[0,T_2]}\left(\mathbf{N}\left(u^{\epsilon,u_0}(t)\right)
-\mathbf{N}\left(u_0\right)\right)\geq\rho\right)\\
\leq&l\exp\left(-(\underline{e}_N-\delta)/\epsilon\right).
\end{array}
\end{equation*}
Thus if $u_0$ belongs to $B_{\rho}^0$
\begin{equation*}
\begin{array}{rl}
\mathbb{P}\left(u^{\epsilon,u_0}\left(\tau^{\epsilon,u_0}\right)\in
N\right)
\leq&\mathbb{P}\left(\tau^{\epsilon,u_0}>\tau_l\right)+\sum_{k=1}^l\mathbb{P}
\left(u^{\epsilon,u_0}\left(\tau^{\epsilon,u_0}\right)\in N,\tau^{\epsilon,u_0}=\tau_k\right)\\
\leq&\mathbb{P}\left(\tau^{\epsilon,u_0}>lT_2\right)+\mathbb{P}\left(\tau_l\leq lT_2\right)\\
&+l\sup_{u_0\in S_{2\rho}^0}\mathbb{P}\left(u^{\epsilon,u_0}
\left(\sigma_{\rho}^{\epsilon,u_0}\right)\in N \right)\\
\leq&\mathbb{P}\left(\tau^{\epsilon,u_0}>lT_2\right)
+3l\exp\left(-(\underline{e}_N-\delta)/\epsilon\right).
\end{array}
\end{equation*}
Take now
$l=\left\lceil(1/T_2)\exp\left((\overline{e}+\delta)/\epsilon\right)\right\rceil$
and use the upper bound \eqref{e3}, possibly choosing $\epsilon_0$
smaller, we obtain that for $\epsilon\leq \epsilon_0$
\begin{equation*}
\begin{array}{rl}
\sup_{u_0\in
B_{\rho}^0}\mathbb{P}\left(u^{\epsilon,u_0}\left(\tau^{\epsilon,u_0}\right)
\in N\right)\leq&\exp\left(-\delta/(5\epsilon)\right)+(4/T_2)
\exp\left(-(\underline{e}_N-\overline{e}+2\delta)/\epsilon\right)\\
\leq&\exp\left(-\delta/(5\epsilon)\right)+(4/T_2)
\exp\left(-\delta/\epsilon\right).
\end{array}
\end{equation*}
Finally, when $u_0$ is any function in $D$, we conclude thanks to
\begin{equation*}
\mathbb{P}\left(u^{\epsilon,u_0}\left(\tau^{\epsilon,u_0}\right)\in
N\right) \leq
\mathbb{P}\left(u^{\epsilon,u_0}\left(\sigma^{\epsilon,u_0}_{\rho}\right)
\in\partial D\right)+ \sup_{u_0\in
B_{\rho}^0}\mathbb{P}\left(u^{\epsilon,u_0}
\left(\tau^{\epsilon,u_0}\right)\in N\right)
\end{equation*}
and to Lemma \ref{l2}.\hfill$\square$
\begin{rmrk}
It is proposed in \cite{SD} to introduce control elements to reduce
or enhance exponentially the expected exit time or to act on the
exiting points, for a limited cost.\@ We could then optimize on
these external fields.\@ However, the problem is computationally
involved since the optimal control problem requires double
optimisation.\@
\end{rmrk}

\section{Exit from a domain of attraction in $\xHone$}\label{s4}
\subsection{Preliminaries}
We now consider a measurable bounded subset $D$ of $\xHone$
invariant by the flow of the deterministic equation and which
contains zero in its interior. We choose $R$ such that $D\subset
B_R^1$.\@ We consider both \eqref{ec1} and \eqref{ec1m} where the
noise is either of additive or of multiplicative type.\@ In this
section we are interested in both the fluctuation of the $\xLtwo$
norm and that of the $\xLtwo$ norm of the gradient.\@ The
Hamiltonian and a modified Hamiltonian are thus of particular
interest.\@ We first distinguish the case where the nonlinearity is
defocusing ($\lambda=-1$) where the Hamiltonian takes non negative
values from the case where the nonlinearity is focusing
($\lambda=1$) where the Hamiltonian may take negative values.\\
\indent We may prove, see for example \cite{GHI}, that
\begin{equation*}
\frac{\xdif}{\xdif t}
\mathbf{H}\left(\mathbf{S}(u_0,0)(t)\right)+2\alpha\mathbf{\Psi}\left(\mathbf{S}(u_0,0)\right)=0,
\end{equation*}
where $\mathbf{S}(u_0,0)$ is the solution of the deterministic
weakly damped nonlinear Schr\"odinger equation with initial datum
$u_0$ in $\xHone$ and
\begin{equation*}
\mathbf{\Psi}\left(\mathbf{S}(u_0,0)\right)=\left\|\nabla
\mathbf{S}(u_0,0)\right\|_{\xLtwo}^2/2-\lambda\int_{\xR^d}
\left|\mathbf{S}(u_0,0)(x)\right|^{2\sigma+2}dx/2.
\end{equation*}
Thus, when the nonlinearity is defocusing we have
\begin{equation}\label{e6}
0\leq\mathbf{H}\left(\mathbf{S}(u_0,0)(t)\right)\leq
\mathbf{H}\left(u_0\right)\exp\left(-2\alpha t\right).
\end{equation}
\indent As it is done in \cite{DO}, we consider in the focusing case
a modified Hamiltonian denoted by $\tilde{\mathbf{H}}(u)$ defined
for $u$ in $\xHone$ by
\begin{equation*}
\tilde{\mathbf{H}}(u)=\mathbf{H}(u)+\beta(\sigma,d)C\left\|u\right\|_{\xLtwo}^{2+4\sigma/(2-\sigma
d)}
\end{equation*}
where the constant $C$ is that of the third inequality in the
following sequence of inequalities where we use the
Gagliardo-Nirenberg inequality
\begin{equation*}
\|u\|_{\xLn^{2\sigma+2}}^{2\sigma+2}/(2\sigma+2)\leq C
\|u\|_{\xLtwo}^{2\sigma+2-\sigma d}\|\nabla u\|_{\xLtwo}^{\sigma
d}\leq\|\nabla u\|_{\xLtwo}^2/4+C\|u\|_{\xLtwo}^{2+4\sigma/(2-\sigma
d)},
\end{equation*}
and $\beta(\sigma,d)=\frac{2\sigma(2-\sigma d)}{(\sigma+2)(2-\sigma
d)+2\sigma(4\sigma+3)}\vee 2$.\@ When evaluated at the deterministic
solution, the modified Hamiltonian satisfies
\begin{equation}\label{e7}
0\leq\tilde{\mathbf{H}}\left(\mathbf{S}(u_0,0)(t)\right)
\leq\tilde{\mathbf{H}}\left(u_0\right)\exp\left(-2\alpha\frac{3(\sigma+1)}{4\sigma+3}
t\right).
\end{equation}
Also, when the nonlinearity is defocusing we now have, for every
$\beta$ positive,
\begin{equation}\label{e8}
0\leq\tilde{\mathbf{H}}\left(\mathbf{S}(u_0,0)(t)\right)
\leq\tilde{\mathbf{H}}\left(u_0\right)\exp\left(-2\alpha t\right).
\end{equation}
From the Sobolev inequalities, for $\rho$ positive, the sets
\begin{equation*}
\tilde{\mathbf{H}}_{\rho}=\left\{u\in\xHone:\
\tilde{\mathbf{H}}(u)=\rho\right\}=\tilde{\mathbf{H}}^{-1}\left(\{\rho\}\right),\quad
\rho>0
\end{equation*}
are closed subsets of $\xHone$ and
\begin{equation*}
\tilde{\mathbf{H}}_{<\rho}=\left\{u\in\xHone:\
\tilde{\mathbf{H}}(u)<\rho\right\}=\tilde{\mathbf{H}}^{-1}\left([0,\rho)\right)\quad
\rho>0
\end{equation*}
are open subsets of $\xHone$.\\
\indent Also, $\tilde{\mathbf{H}}$ is such that
\begin{equation}\label{e9}
\|\nabla u\|_{\xLtwo}^2/2+\beta C\|u\|_{\xLtwo}^{2+4\sigma/(2-\sigma
d)} \leq\tilde{\mathbf{H}}(u)\leq3\|\nabla u\|_{\xLtwo}^2/4+(\beta
+1)C\|u\|_{\xLtwo}^{2+4\sigma/(2-\sigma d)}
\end{equation}
when the nonlinearity is defocusing and
\begin{equation}\label{e10}
\|\nabla u\|_{\xLtwo}^2/4+ C\|u\|_{\xLtwo}^{2+4\sigma/(2-\sigma d)}
\leq\tilde{\mathbf{H}}(u)\leq\|\nabla
u\|_{\xLtwo}^2/2+\beta(\sigma,d)C\|u\|_{\xLtwo}^{2+4\sigma/(2-\sigma
d)}
\end{equation}
when it is focusing.\@ Thus the sets $\tilde{\mathbf{H}}_{<\rho}$
for $\rho$ positive are bounded in $\xHone$ and a bounded set in
$\xHone$ is bounded for $\tilde{\mathbf{H}}$.\@ Note that the domain
$D$ of attraction may be a domain of the form
$\tilde{\mathbf{H}}_{<\rho}$.\@ \vspace{0.3cm}

\indent We no longer distinguish the focusing and defocusing cases
and take the same value of $\beta$, {\it i.e.} $\beta(\sigma,d)$.\@
Also to simplify the notations we now sometimes drop the dependence
of the solution
in $\epsilon$ and $u_0$.\\
\indent The fluctuation of
$\tilde{\mathbf{H}}\left(u^{\epsilon,u_0}(t)\right)$ is of
particular interest.\@ We have the following result when the noise
is of additive type.\@
\begin{prpstn}
When $u$ denotes the solution of equation \eqref{ec1},
$\left(e_j\right)_{j\in\xN}$ a complete orthonormal system of
$\xLtwo$, the following decomposition holds
\begin{equation*}
\begin{array}{rl}
\tilde{\mathbf{H}}\left(u(t)\right)=&\tilde{\mathbf{H}}\left(u_0\right)\\
&-2\alpha\int_0^t\mathbf{\Psi}\left(u(s)\right)ds-2\beta
C\left(1+2\sigma/(2-\sigma
d)\right)\alpha\int_0^t\|u(s)\|_{\xLtwo}^{2+4\sigma/(2-\sigma
d)}ds\\
&+\sqrt{\epsilon}\left(\pim\int_{\xR^d}\int_0^t\nabla\overline{u}(s)\nabla
dW(s)dx-\lambda\pim\int_{\xR^d}\int_0^t\left|u(s)\right|^{2\sigma}\overline{u}(s)dW(s)dx\right.\\
&\left.\quad\quad\quad+2\beta C\left(1+2\sigma/(2-\sigma
d)\right)\pim\int_{\xR^d}\int_0^t\left\|u(s)\right\|_{\xLtwo}^{4\sigma/(2-\sigma
d)}\overline{u}(s)dW(s)dx\right)\\
&-(\lambda\epsilon/2)\sum_{j\in\xN}\int_0^t\int_{\xR^d}\left[\left|u(s)\right|^{2\sigma}|\Phi
e_j|^2+2\sigma\left|u(s)\right|^{2\sigma-2}(\pre(\overline{u}(s)\Phi
e_j))^2\right]dxds\\
&+(\epsilon/2)\|\nabla\Phi\|_{\mathcal{L}_2^{0,0}}^2 t+\epsilon\beta
C\left(1+2\sigma/(2-\sigma
d)\right)\|\Phi\|_{\mathcal{L}_2^{0,0}}^2\int_0^t\|u(s)\|_{\xLtwo}^{4\sigma/(2-\sigma
d)}ds\\
&+\epsilon\beta C\left(4\sigma/(2-\sigma
d)\right)\left(1+2\sigma/(2-\sigma
d)\right)\sum_{j\in\xN}\int_0^t\|u(s)\|_{\xLtwo}^{2\left(2\sigma/(2-\sigma
d)-1\right)}\left(\pre\int_{\xR^d}\overline{u}(s)\Phi
e_jdx\right)^2ds
\end{array}
\end{equation*}
\end{prpstn}
{\bf Proof.}  The result follows from the It\^o formula.\@ The main
difficulty is in justifying the computations.\@ We may proceed as in
\cite{dBD1}.\hfill $\square$\vspace{0.3cm}

\indent Also, when the noise is of multiplicative type we obtain the
following proposition.\@
\begin{prpstn}
When $u$ denotes the solution of equation \eqref{ec1m},
$\left(e_j\right)_{j\in\xN}$ a complete orthonormal system of
$\xLtwo$, the following decomposition holds
\begin{equation*}
\begin{array}{ll}
\tilde{\mathbf{H}}\left(u(t)\right)=&\tilde{\mathbf{H}}\left(u_0\right)\\
&-2\alpha\int_0^t\mathbf{\Psi}\left(u(s)\right)ds-2\beta
C\left(1+2\sigma/(2-\sigma
d)\right)\alpha\int_0^t\|u(s)\|_{\xLtwo}^{2+4\sigma/(2-\sigma
d)}ds\\
&+\sqrt{\epsilon}\pim\int_{\xR^d}\int_0^tu(s)\nabla\overline{u}(s)\nabla
dW(s)dx\\
&+(\epsilon/2)\sum_{j\in\xN}\int_0^t\int_{\xR^d}|u(s)|^2|\nabla\Phi
e_j|^2dxds.
\end{array}
\end{equation*}
\end{prpstn}
\indent The first exit time $\tau^{\epsilon,u_0}$ from the domain
$D$ in $\xHone$ is defined as in Section \ref{s2}.\@ We also define
\begin{equation*}
\overline{e}=\inf\left\{I_T^{0}(w):\ w(T)\in \overline{D}^c,\
T>0\right\},
\end{equation*}
and for $\rho$ positive small enough
\begin{equation*}
e_{\rho}=\inf\left\{I_T^{u_0}(w):\
\tilde{\mathbf{H}}\left(u_0\right)\leq \rho,\ w(T)\in
\left(D_{-\rho}\right)^c,\ T>0\right\},
\end{equation*}
where $D_{-\rho}=D\setminus\mathcal{N}^1\left(\partial
D,\rho\right)$.\@ Then we set
\begin{equation*}
\underline{e}=\lim_{\rho\rightarrow0}e_{\rho}.\end{equation*} Also,
for $\rho$ positive small enough, $N$ a closed subset of the
boundary of $D$, we define
\begin{equation*}
e_{N,\rho}=\inf\left\{I_T^{u_0}(w):\
\tilde{\mathbf{H}}\left(u_0\right)\leq\rho,\ w(T)\in
\left(D\setminus \mathcal{N}^1\left(N,\rho\right)\right)^c,\
T>0\right\}
\end{equation*}
and
\begin{equation*}
\underline{e}_N=\lim_{\rho\rightarrow0}e_{N,\rho}.
\end{equation*}
We finally also introduce
\begin{equation*}
\sigma_{\rho}^{\epsilon,u_0}=\inf\left\{t\geq
0:u^{\epsilon,u_0}(t)\in \tilde{\mathbf{H}}_{<\rho}\cup D^c\right\},
\end{equation*}
where $\tilde{\mathbf{H}}_{<\rho}\subset D$.\\
\indent Again we have the following inequalities.\@
\begin{lmm}\label{l02}
$0<\underline{e}\leq\overline{e}$.\@
\end{lmm}
{\bf Proof.}  We only have to prove the first inequality.\@
Integrating the equation describing the evolution of
$\tilde{\mathbf{H}}\left(\mathbf{S}\left(u_0,h\right)(t)\right)$ via
the Duhamel formula where the skeleton is that of the equation with
an additive noise we obtain
\begin{equation*}
\begin{array}{l}
\tilde{\mathbf{H}}\left(\mathbf{S}(u_0,h)(T)\right)-
\exp\left(-2\alpha\frac{3(\sigma+1)}{4\sigma+3}
T\right)\tilde{\mathbf{H}}\left(u_0\right)\\
\leq\int_0^T\exp\left(-2\alpha\frac{3(\sigma+1)}{4\sigma+3}
(T-s)\right)\left[\pim\int_{\xR^d}\left(\nabla\mathbf{S}(u_0,h)\nabla\overline{\Phi
h}\right)(s,x)dx\right.\\
\left.\quad
-\lambda\pim\int_{\xR^d}\left(|\mathbf{S}(u_0,h)|^{2\sigma}\mathbf{S}(u_0,h)\overline{\Phi
h}\right)(s,x)dx\right.\\
\left.\quad -2C\beta\left(1+2\sigma/(2-\sigma
d)\right)\pim\int_{\xR^d}\left(\mathbf{S}\left(u_0,h\right)\overline{\Phi
h}\right)(s,x)dx\right]ds,
\end{array}
\end{equation*}
with a focusing or defocusing nonlinearity.\@ Let $d$ denote the
positive distance between 0 and $\partial D$.\@ Take $\rho$ such
that the distance between $B_{\rho}^1$ and
$\left(D_{-\rho}\right)^c$ is larger than $d/2$.\@ We then have,
from the fact that the Sobolev injection from $\xHone$ into
$\xLn^{2\sigma+2}$,
\begin{equation*}
\begin{array}{l}
d/2\leq\int_0^T\exp\left(-2\alpha\frac{3(\sigma+1)}{4\sigma+3}
(T-s)\right)\left[R\|\Phi\|_{\mathcal{L}_c\left(\xLtwo,\xHone\right)}\|h\|_{\xLtwo}\right.\\
\left.\quad
+CR^{2\sigma+1}\|\Phi\|_{\mathcal{L}_c\left(\xLtwo,\xHone\right)}\|h\|_{\xLtwo}\right.\\
\left.\quad +2C\beta\left(1+2\sigma/(2-\sigma
d)\right)R\|\Phi\|_{\mathcal{L}_c\left(\xLtwo,\xLtwo\right)}\|h\|_{\xLtwo}\right]ds,
\end{array}
\end{equation*}
We conclude as in Lemma \ref{l0} and use that from the choice of
$\beta$ the complementary of a ball is included in the complementary
of a set $\tilde{\mathbf{H}}_{<a}$.\@ In the case of the skeleton of
the equation with a multiplicative noise, it is enough to replace
the term in bracket in the right hand side of the above formula by
$\pim\int_{\xR^d}\left(\nabla\mathbf{S}\left(u_0,h\right)\overline{\mathbf{S}\left(u_0,h\right)}\overline{\nabla\Phi
h}\right)(s,x)dx$.\@ Recall that we can proceed as in the additive
case since we have imposed that $\Phi$ belongs to
$\mathcal{L}_{2,\xR}^{0,s}$ where $s>d/2+1$, in particular $\Phi$
belongs to $\mathcal{L}_c\left(\xLtwo,\xW^{1,\infty}\right)$.\hfill
$\square$

\subsection{Statement of the results}
The theorems of Section \ref{s2} still hold for a domain of
attraction in $\xHone$ and a noise of additive and multiplicative
type.\@
\begin{thm}\label{t22}
For every $u_0$ in $D$ and $\delta$ positive, there exists $L$
positive such that
\begin{equation}\label{et212}
\overline{\lim}_{\epsilon\rightarrow0}\epsilon\log\mathbb{P}
\left(\tau^{\epsilon,u_0}\notin\left(\exp\left((\underline{e}
-\delta)/\epsilon\right),\exp\left((\overline{e}
+\delta)/\epsilon\right)\right)\right)\leq-L,
\end{equation}
and for every $u_0$ in $D$,
\begin{equation}\label{et222}
\underline{e}\leq\underline{\lim}_{\epsilon\rightarrow0}\epsilon\log
\mathbb{E}\left(\tau^{\epsilon,u_0}\right)
\leq\overline{\lim}_{\epsilon\rightarrow0}\epsilon\log
\mathbb{E}\left(\tau^{\epsilon,u_0}\right)\leq\overline{e}.
\end{equation}
Moreover, for every $\delta$ positive, there exists $L$ positive
such that
\begin{equation}\label{et232}
\overline{\lim}_{\epsilon\rightarrow0}\epsilon\log\sup_{u_0\in
D}\mathbb{P}\left(\tau^{\epsilon,u_0}\geq
\exp\left((\overline{e}+\delta)/\epsilon\right)\right)\leq-L,
\end{equation}
and
\begin{equation}\label{et242}
\overline{\lim}_{\epsilon\rightarrow0}\epsilon\log\sup_{u_0\in D}
\mathbb{E}\left(\tau^{\epsilon,u_0}\right)\leq\overline{e}.
\end{equation}
\end{thm}
\begin{rmrk}
Again the control argument to prove that
$\underline{e}=\overline{e}$ seems difficult.\@ It should be even
more difficult for multiplicative noises.
\end{rmrk}
\begin{thm}\label{t32} If
$\underline{e}_N>\overline{e}$, then for every $u_0$ in $D$, there
exists $L$ positive such that
\begin{equation*}
\overline{\lim}_{\epsilon\rightarrow
0}\epsilon\log\mathbb{P}\left(u^{\epsilon,u_0}\left(\tau^{\epsilon,u_0}\right)\in
N\right)\leq-L.
\end{equation*}
\end{thm}
Again we may deduce the corollary
\begin{crllr}
Assume that $v^*$ in $\partial D$ is such that for every $\delta$
positive and $N=\left\{v\in\partial D:\ \|v-v^*\|_{\xLtwo}\geq
\delta\right\}$ we have $\underline{e}_N>\overline{e}$ then
\begin{equation*}
\forall \delta>0,\ \forall u_0\in D,\ \exists L>0:\
\overline{\lim}_{\epsilon\rightarrow0}\epsilon
\log\mathbb{P}\left(\left\|u^{\epsilon,u_0}
\left(\tau^{\epsilon,u_0}\right)-v^*\right\|_{\xLtwo}\geq\delta\right)\leq-L.
\end{equation*}
\end{crllr}
\subsection{Proof of the results}
The proof of these results still rely on three lemmas and the
uniform LDP.\@ Let us now state the lemmas for both a noise of
additive and of multiplicative type.\@
\begin{lmm}\label{l12} For every $\rho$ and $L$ positive with
$\tilde{\mathbf{H}}_{<\rho}\subset D$, there exists $T$ and
$\epsilon_0$ positive such that for every $u_0$ in $D$ and
$\epsilon$ in $(0,\epsilon_0)$,
\begin{equation*}
\mathbb{P}\left(\sigma_{\rho}^{\epsilon,u_0}>T\right)\leq\exp\left(-L/\epsilon\right).
\end{equation*}
\end{lmm}
{\bf Proof.}  We proceed as in the proof of Lemma \ref{l1}.\\
Let $d$ denote the positive distance between 0 and $D\setminus
\tilde{\mathbf{H}}_{<\rho}$.\@ Take $\alpha$ positive such that
$\alpha\rho<d$.\@ The domain $D$ is uniformly attracted to 0, thus
there exists a time $T_1$ such that for every initial datum $u_1$ in
$\mathcal{N}^1\left(D\setminus
\tilde{\mathbf{H}}_{<\rho},\alpha\rho/8\right)$, for $t\geq T_1$,
$\mathbf{S}\left(u_1,0\right)(t)$ belongs to
$B_{\alpha\rho/8}^1$.\\
\indent We could also prove, see \cite{dBD1}, that there exists a
constant $M'$ which depends on $T_1$, $R$, $\sigma$ and $\alpha$
such that
\begin{equation}\sup_{u_1\in\mathcal{N}^1\left(D\setminus
\tilde{\mathbf{H}}_{<\rho},\alpha\rho/8\right)}
\left\|\mathbf{S}\left(u_1,0\right)\right\|_{X^{(T_1,2\sigma+2)}}\leq
M'.\end{equation}\indent The Step 2, corresponding to that of Lemma
\ref{l1}, in the proof in the additive case uses the truncation
argument, upper bounds similar to that in \cite{dBD1} derived from
the Strichartz inequalities on smaller intervals; we shall also
replace in the proof of Lemma \ref{l1} $\rho/8$ by
$\alpha\rho/8$.\\
\indent In Step 2 for the multiplicative case, we also introduce the
truncation in front of the term $u\Phi h$ in the controlled PDE.\\
\indent The end of the proof is identical to that of Lemma \ref{l1},
the LDP is the LDP in $\xCn\left([0,T];\xHone\right)$, for additive
or multiplicative noises.\hfill $\square$
\begin{lmm}\label{l22}
For every $\rho$ positive such that
$\tilde{\mathbf{H}}_{\rho}\subset D$ and $u_0$ in $D$, there exists
$L$ positive such that
\begin{equation*}
\overline{\lim}_{\epsilon\rightarrow 0}\epsilon\log
\mathbb{P}\left(u^{\epsilon,u_0}\left(\sigma_{\rho}^{\epsilon,u_0}\right)\in\partial
D \right)\leq-L
\end{equation*}
\end{lmm}
{\bf Proof.} It is the same proof as for Lemma \ref{l2}.\@ We only
have to replace $B_{\rho/2}^0$ by any ball in $\xHone$ centered at 0
and included in $\tilde{\mathbf{H}}_{<\rho}$ and use the LDP in
$\xCn\left([0,T];\xHone\right).\square$
\begin{lmm}\label{l32}
For every $\rho$ and $L$ positive such that
$\tilde{\mathbf{H}}_{2\rho}\subset D$, there exists
$T(L,\rho)<\infty$ such that
\begin{equation*}
\overline{\lim}_{\epsilon\rightarrow0}\epsilon\log\sup_{u_0\in
\tilde{\mathbf{H}}_{\rho}}\mathbb{P}\left(\sup_{t\in[0,T(L,\rho)]}
\left(\tilde{\mathbf{H}}\left(u^{\epsilon,u_0}(t)\right)
-\tilde{\mathbf{H}}\left(u_0\right)\right)\geq \rho\right)\leq-L
\end{equation*}
\end{lmm}
{\bf Proof.} Integrating the It\^o differential relation using the
Duhamel formula allows to get rid of the drift term that is not
originated from the bracket.\@ Indeed, the event
\begin{equation*}
\left\{\sup_{t\in[0,T(L,\rho)]}\left(\tilde{\mathbf{H}}
\left(u^{\epsilon,u_0}(t)\right)-\tilde{\mathbf{H}}\left(u_0\right)\right)\geq
\rho\right\}
\end{equation*}
is included in
\begin{equation*}
\left\{\sup_{t\in[0,T(L,\rho)]}\left(\tilde{\mathbf{H}}
\left(u^{\epsilon,u_0}(t)\right)-\exp\left(-2\alpha
\left(\frac{3(\sigma+1)}{4\sigma+3}\right)T(L,\rho)\right)
\tilde{\mathbf{H}}\left(u_0\right)\right)\geq \rho\right\}.
\end{equation*}
Then, setting $c(\sigma)=\frac{3(\sigma+1)}{4\sigma+3}$ and
$m(\sigma,d)=1+2\sigma/(2-\sigma d)$, dropping the exponents
$\epsilon$ and $u_0$ to have more concise formulas, we obtain in the
additive case
\begin{equation*}
\begin{array}{l}
\tilde{\mathbf{H}}\left(u(t)\right)-\exp\left(-2\alpha
c(\sigma)t\right)
\tilde{\mathbf{H}}\left(u_0\right)\\
\leq\sqrt{\epsilon}\left(\pim\int_{\xR^d}\int_0^t\exp\left(-2\alpha
c(\sigma)(t-s)\right)\nabla\overline{u}(s)\nabla
dW(s)dx\right.\\
\left.-\lambda\pim\int_{\xR^d}\int_0^t\exp\left(-2\alpha
c(\sigma)(t-s)\right)\left|u(s)\right|^{2\sigma}\overline{u}(s)dW(s)dx\right.\\
\left.+2\beta Cm(\sigma,d)\pim\int_{\xR^d}\int_0^t\exp\left(-2\alpha
c(\sigma)(t-s)\right)\left\|u(s)\right\|_{\xLtwo}^{4\sigma/(2-\sigma
d)}\overline{u}(s)dW(s)dx\right)\\
-(\lambda\epsilon/2)\sum_{j\in\xN}\int_0^t\exp\left(-2\alpha
c(\sigma)(t-s)\right)\int_{\xR^d}\left[\left|u(s)\right|^{2\sigma}|\Phi
e_j|^2\right.\\
\quad\quad\quad\quad\quad\quad\quad\quad\quad\quad\quad\quad\quad\quad\quad\quad\quad+\left.2\sigma\left|u(s)\right|^{2\sigma-2}(\pre(\overline{u}(s)\Phi
e_j))^2\right]dxds\\
+\left(\epsilon/(4\alpha c(\sigma))\right)\left(1-\exp\left(-2\alpha
c(\sigma)t\right)\right)\|\nabla\Phi\|_{\mathcal{L}_2^{0,0}}^2\\
+\epsilon\beta
Cm(\sigma,d)\|\Phi\|_{\mathcal{L}_2^{0,0}}^2\int_0^t\exp\left(-2\alpha
c(\sigma)(t-s)\right)\|u(s)\|_{\xLtwo}^{4\sigma/(2-\sigma
d)}ds\\
+\epsilon\beta Cm(\sigma,d)\left(4\sigma/(2-\sigma
d)\right)\sum_{j\in\xN}\int_0^t\exp\left(-2\alpha
c(\sigma)(t-s)\right)\|u(s)\|_{\xLtwo}^{2\left(2\sigma/(2-\sigma
d)-1\right)}\left(\pre\int_{\xR^d}\overline{u}(s)\Phi
e_jdx\right)^2ds.
\end{array}
\end{equation*}
We again use a localization argument and replace the process $u$ by
the process $u^{\tau}$ stopped at the first exit time off
$\tilde{\mathbf{H}}_{<2\rho}$.\@ We use \eqref{e9} and \eqref{e10}
and obtain
\begin{equation*}
\|u^{\tau}\|_{\xHone}^2\leq 8\rho+\left(2\rho/(C\sigma)
\right)^{\frac{1}{1+2\sigma/(2-\sigma d)}}.
\end{equation*}
We denote the right hand side of the above by $b(\rho,\sigma,d)$.\\
From the H\"older inequality along with the Sobolev injection of
$\xHone$ into $\xLn^{2\sigma+2}$ we obtain the following upper bound
for the drift
\begin{equation*}
\begin{array}{l}
\left(\epsilon/(4\alpha
c(\sigma))\right)\left[(1+2\sigma)c(1,2\sigma+2)^{2\sigma+2}\|\Phi\|_{\mathcal{L}_2^{0,1}}^2b(\rho,\sigma,d)^{2\sigma}+\|\nabla\Phi\|_{\mathcal{L}_2^{0,0}}^2\right]\\
+m(\sigma,d)\left(\epsilon\beta C/(2\alpha
c(\sigma))\right)\left(1+4\sigma/(2-\sigma
d)\right)\|\Phi\|_{\mathcal{L}_2^{0,0}}^2b(\rho,\sigma,d)^{4\sigma/(2-\sigma
d)}
\end{array}
\end{equation*}
where we denote by $c(1,2\sigma+2)$ the norm of the continuous
injection of $\xHone$ into $\xLn^{2\sigma+2}$.\\
Thus, choosing $\epsilon$ small enough, it is enough to show the
result for the stochastic integral replacing $\rho$ by $\rho/2$.\@
Also it is enough to show the result for each of the three
stochastic integrals replacing $\rho/2$ by $\rho/6$.\@ With the same
one parameter families and similar computations as in the proof of
Lemma \ref{l3}, we know that it is enough to obtain upper bounds of
the brackets of the stochastic integrals
\begin{equation*}
\begin{array}{l}
Z_1(t)=\pim\int_{\xR^d}\int_0^t\exp\left(2\alpha
c(\sigma)s\right)\nabla\overline{u^{\tau}}(s)\nabla
dW(s)dx\\
Z_2(t)=\pim\int_{\xR^d}\int_0^t\exp\left(2\alpha c(\sigma)s\right)
\left|u^{\tau}(s)\right|^{2\sigma}\overline{u^{\tau}}(s)dW(s)dx\\
Z_3(t)=2\beta Cm(\sigma,d)\pim\int_{\xR^d}\int_0^t\exp\left(2\alpha
c(\sigma)s\right)\left\|u^{\tau}(s)\right\|_{\xLtwo}^{4\sigma/(2-\sigma
d)}\overline{u^{\tau}}(s)dW(s)dx.
\end{array}
\end{equation*}
We then obtain
\begin{equation*}
\begin{array}{l}
d<Z_1>_t\leq\exp\left(4\alpha c(\sigma)t\right)\sum_{j\in\xN}
\left(\nabla u^{\tau}(t),-i\nabla\Phi e_j\right)_{\xLtwo}^2dt\\
d<Z_2>_t\leq\exp\left(4\alpha c(\sigma)t\right)\sum_{j\in\xN}
\left(|u^{\tau}(t)|^{2\sigma}u^{\tau}(t),-i\Phi e_j\right)_{\xLtwo}^2dt\\
d<Z_3>_t\leq 4\beta^2 C^2m(\sigma,d)^2\exp\left(4\alpha
c(\sigma)t\right)\left\|u^{\tau}(t)\right\|_{\xLtwo}^{8\sigma/(2-\sigma
d)}\sum_{j\in\xN}\left(u^{\tau}(t),-i\Phi e_j\right)_{\xLtwo}^2dt.
\end{array}
\end{equation*}
Using the H\"older inequality and, for $Z_2$, the continuous Sobolev
injection of $\xHone$ into $\xLn^{2\sigma+2}$ we obtain
\begin{equation*}
\begin{array}{l}
d<Z_1>_t\leq\exp\left(4\alpha c(\sigma)t\right)\|\Phi\|_{\mathcal{L}_2^{0,1}}^2b(\rho,\sigma,d)dt\\
d<Z_2>_t\leq\exp\left(4\alpha
c(\sigma)t\right)c(1,2\sigma+2)^{2(2\sigma+2)}\|\Phi\|_{\mathcal{L}_2^{0,1}}^2
b(\rho,\sigma,d)^{2\sigma+1}dt\\
d<Z_3>_t\leq4\beta^2 C^2m(\sigma,d)^2\exp\left(4\alpha
c(\sigma)t\right)b(\rho,\sigma,d)^{\left(1+4\sigma/(2-\sigma
d)\right)}\|\Phi\|_{\mathcal{L}_2^{0,1}}^2dt.
\end{array}
\end{equation*}
We can then bound each of the three remainders
$\left(R_l^i(t)\right)_{i=1,2,3}$ similar to that of Lemma \ref{l3}
using the inequality
$R_l^i(t)\leq3l\int_0^td<Z_i>_t$.\\
We conclude that it is possible to choose $T(L,\rho)$ equal to\\
$\frac{1}{4\alpha c(\sigma)}\log\left(\frac{\alpha c(\sigma)
\rho^2}{90 b(\rho,\sigma,d)\|\Phi\|_{\mathcal{L}_2^{0,1}}^2
\max\left(1,c(1,2\sigma+2)^{2(2\sigma+1)}b(\rho,\sigma,d)^{2\sigma},4\beta^2
C^2m(\sigma,d)^2b(\rho,\sigma,d)^{4\sigma/(2-\sigma d)}\right)}\right).$\\
When the noise is of multiplicative type we obtain
\begin{equation*}
\begin{array}{l}
\tilde{\mathbf{H}}\left(u(t)\right)-\exp\left(-2\alpha
c(\sigma)t\right)\tilde{\mathbf{H}}\left(u_0\right)\\
\leq\sqrt{\epsilon}\pim\int_{\xR^d}\int_0^t\exp\left(-2\alpha
c(\sigma)(t-s)\right)u(s)\nabla\overline{u}(s)\nabla
dW(s)dx\\
\quad+(\epsilon/2)\sum_{j\in\xN}\int_0^t\exp\left(-2\alpha
c(\sigma)(t-s)\right)\int_{\xR^d}|u(s)|^2|\nabla\Phi e_j|^2dxds.
\end{array}
\end{equation*}
Again we use a localization argument and consider the process $u$
stopped at the exit off $\tilde{\mathbf{H}}_{2\rho}$.\@ As $\Phi$ is
Hilbert-Schmidt from $\xLtwo$ into $\xHn_{\xR}^s$, the second term
of the right hand side is less than $\frac{\epsilon}{4\alpha
c(\sigma)}\|\Phi\|_{\mathcal{L}_2^{0,s}}^2b(\rho,\sigma,d)$ and for
$\epsilon$ small enough, it is enough to prove the result for the
stochastic integral replacing $\rho$ by $\rho/2$.\@ We know that it
is enough to obtain an upper bound of the bracket of
\begin{equation*}
Z(t)=\pim\int_{\xR^d}\int_0^t\exp\left(2\alpha c(\sigma)s\right)
u^{\tau}(s)\nabla\overline{u}^{\tau}(s)\nabla dW(s)dx.
\end{equation*}
We obtain
\begin{equation*}
d<Z>_t\leq\exp\left(4\alpha
c(\sigma)t\right)\sum_{j\in\xN}\left(\nabla u^{\tau}(t),
-iu^{\tau}(t)\nabla\Phi e_j\right)_{\xLtwo}^2dt.
\end{equation*}
Denoting by $c(s,\infty)$ the norm of the Sobolev injection of
$\xHn^s_{\xR}$ into $\xW^{1,\infty}_{\xR}$ we deduce that
\begin{equation*}
d<Z>_t\leq\exp\left(4\alpha c(\sigma)t\right)
c(s,\infty)^2\|\Phi\|_{\mathcal{L}_2^{0,s}}^2b(\rho,\sigma,d)^2dt.
\end{equation*}
Finally, we conclude that we may choose
\begin{equation*}T(L,\rho)=\frac{1}{4\alpha c(\sigma)}\log\left(\frac{\alpha
c(\sigma)\rho^2}{10
b(\rho,\sigma,d)^2c(s,\infty)^2\|\Phi\|_{\mathcal{L}_2^{0,s}}^2L}\right).\end{equation*}\hfill
$\square$\vspace{0.3cm}

\indent We may now prove Theorem \ref{et212} and \ref{et222}.\\
{\bf Elements for the proof of Theorem \ref{et212}.} There is no
difference in the proof of the upper bound on
$\tau^{\epsilon,u_0}$.\@ Let us thus focus on the lower bound.\@
Take $\delta$ positive.\@ Since $\underline{e}>0$, we now choose
$\rho$ positive such that $\underline{e}-\delta/4\leq e_{\rho}$,
$\tilde{\mathbf{H}}_{2\rho}\subset D$ and
$\tilde{\mathbf{H}}_{2\rho}\subset D_{-\rho}^c$.\@ We define the
sequences of stopping times $\theta_0=0$ and for $k$ in $\xN$,
\begin{equation*}\begin{array}{rl}
\tau_k&=\inf\left\{t\geq\theta_k:\ u^{\epsilon,u_0}(t)\in
\tilde{\mathbf{H}}_{<\rho}\cup
D^c\right\},\\
\theta_{k+1}&=\inf\left\{t>\tau_k:\ u^{\epsilon,u_0}(t)\in
\tilde{\mathbf{H}}_{2\rho}\right\},\end{array}
\end{equation*}
where $\theta_{k+1}=\infty$ if $u^{\epsilon,u_0}(\tau_k)\in
\partial D$.\@ Let us fix $T_1=T\left(\underline{e}-3\delta/4,\rho\right)$
given by Lemma \ref{l32}.\@ We now use that for $u_0$ in $D$ and $m$
a positive integer,
\begin{equation}\label{e42}
\begin{array}{rl}
\mathbb{P}\left(\tau^{\epsilon,u_0}\leq mT_1\right)\leq&
\mathbb{P}\left(\tau^{\epsilon,u_0}=\tau_0\right)+\sum_{k=1}^m
\mathbb{P}\left(\tau^{\epsilon,u_0}=\tau_k\right)\\
&+\sum_{k=1}^m\mathbb{P}\left(\theta_k-\tau_{k-1}\leq T_1\right)
\end{array}\end{equation}
and conclude as in the proof of Theorem \ref{t2}.\hfill
$\square$\vspace{0.3cm}

We may check that the proof of Theorem \ref{t3} also applies to
Theorem \ref{t32}, the LDPs are those in $\xHone$ and the sequences
of stopping times are those defined above.
\begin{rmrk}
In \cite{Fr}, reaction-diffusion equations perturbed by an additive
white noise are considered. When the space dimension is larger than
one, the case where the vector field can decomposed in a gradient
and a second field which is orthogonal is treated. The
quasi-potential is then equal to the potential at the end point. It
again involves a control argument. In our case, since we consider
colored noises and nonlinear equations, the orthogonality is lost
for the geometry of the reproducing kernel Hilbert space of the law
of $W(1)$. We thus obtain extra commutator terms. Under suitable
assumptions on the space correlations of the noise, going to zero,
it is possible that we obtain a non trivial minimisation problem.
Recall that solitary waves are solutions of variational problem
where we minimize the Hamiltonian for fixed levels of the mass.
\end{rmrk}

\section{Annex - proof of Theorem \ref{t1}}\label{s5}
The following lemma is at the core of the proof of the uniform
LDPs.\@ It is often called Azencott lemma or Freidlin-Wentzell
inequality.\@ The differences with the result of \cite{EG2} are that
here the initial data are the same for the random process and the
skeleton and that the "for every $\rho$ positive" stands before
"there exists $\epsilon_0$ and $\gamma$ positive".\@ We shall only
stress on the differences in the proof.\@
\begin{lmm}\label{6l4}
For every $a$, $L$, $T$, $\delta$ and $\rho$ positive, $f$ in $C_a$,
$p$ in $\mathcal{A}(d)$, there exists $\epsilon_0$ and $\gamma$
positive such that for every $\epsilon$ in $(0,\epsilon_0)$,
$\|u_0\|_{\xHone}\leq \rho$,
\begin{equation*}
\epsilon\log\mathbb{P}\left(\left\|u^{\epsilon,u_0}-
\tilde{\mathbf{S}}(u_0,f)\right\|_{X^{(T,p)}} \geq\delta;\
\|\sqrt{\epsilon}W-f\|_{\xCn\left([0,T]; \xHn_{\xR}^s\right)}<\gamma
\right)\leq-L.
\end{equation*}
\end{lmm}
{\bf Elements of proof.} There are still three steps in the proof of
this result.\@ The first step is a change of measure to center the
process
around $f$.\@ It uses the Girsanov theorem and is the same as in \cite{EG2}.\\
The second step is a reduction to estimates for the stochastic
convolution.\@ It strongly involves the Strichartz inequalities but
it is slightly different than in \cite{EG2}.\@ The truncation
argument has to hold for all $\|u_0\|_{\xHone}\leq \rho$.\@ Thus we
use the fact that there exists $M=M(T,\rho,\sigma)$ positive such
that
\begin{equation*}\sup_{u_1\in
B_{\rho}^1}\left\|\tilde{\mathbf{S}}(u_1,f)
\right\|_{X^{\left(T,p\right)}}\leq M.\end{equation*} The proof of
this fact follows from the computations in \cite{dBD1}, we have
recalled the arguments in $\xLtwo$ in the proof of Lemma \ref{l1}.\@
The result in $\xHone$ is again be used in the proof of Lemma
\ref{l12}.\@ As the initial data are the same for the random process
and the skeleton, the remaining of the argument
does not require restrictions on $\rho$.\\
The third step corresponds to estimates for the stochastic
convolution.\@ It is the same as in \cite{EG2}.\\
The extra damping term in the drift is treated easily thanks to the
Strichartz inequalities.\hfill $\square$\vspace{0.3cm}

{\bf Elements for the proof of Theorem \ref{t1}.} Let us start with
the case of an additive noise.\@ Recall that, in that case, the mild
solution of the stochastic equation could be written as a function
of the perturbation in the convolution form.\@ Let $v^{u_0}(Z)$
denote the solution of
\begin{equation*}\left\{\begin{array}{l}
    i\frac{\partial v}{\partial t} -\left(\Delta v+
|v-iZ|^{2\sigma}(v-iZ)-i\alpha(v-iZ)\right)=0, \\
    v(0)=u_0,
\end{array}\right.\end{equation*}
or equivalently a fixed point of the functional $\mathcal{F}_{Z}$
such that
\begin{equation*}
\begin{array}{rl}
\mathcal{F}_{Z}(v)(t)=&U(t)u_0-i\lambda\int_0^tU(t-s)
\left(|(v-iZ)(s)|^{2\sigma}(v-iZ)(s)\right)ds\\
&-\alpha\int_0^tU(t-s)(v-iZ)(s)ds,
\end{array}
\end{equation*}
where $Z$ belongs to $\xCn\left([0,T];\xLtwo\right)$ (respectively
$\xCn\left([0,T];\xHone\right)$).\@ If $u^{\epsilon,u_0}$ is defined
as $u^{\epsilon,u_0}
=v^{u_0}\left(Z^{\epsilon}\right)-iZ^{\epsilon}$ where
$Z^{\epsilon}$ is the stochastic convolution
$Z^{\epsilon}(t)=\sqrt{\epsilon}\int_0^tU(t-s)dW(s)$ then
$u^{\epsilon,u_0}$ is a solution of the stochastic equation.\@
Consequently, if $\mathcal{G}\left(\cdot,u_0\right)$ denotes the
mapping from $\xCn\left([0,T];\xLtwo\right)$ (respectively
$\xCn\left([0,T];\xHone\right)$) to $\xCn\left([0,T];\xLtwo\right)$
(respectively $\xCn\left([0,T];\xHone\right)$) defined by
$\mathcal{G}\left(Z,u_0\right)=v^{u_0}(Z)-iZ$, we obtain
$u^{\epsilon,u_0}=\mathcal{G}\left(Z^{\epsilon},u_0\right)$.\@ We
may also check with arguments similar to that of \cite{dBD1,EG1},
involving the Strichartz inequalities that the mapping $\mathcal{G}$
is equicontinuous in its first arguments for second arguments in
bounded sets of $\xLtwo$ (respectively $\xHone$).\@ The result now
follows from Proposition 5 in \cite{Sow}.\vspace{0.3cm}

\indent Let us now consider the case of a multiplicative noise.\@
Initial data belong to $\xHone$ and we consider paths in $\xHone$.\@
The proof is very
close to that in \cite{EG2}.\\
\indent The main tool is again the Azencott lemma or almost
continuity of the It\^o map.\@ We need the slightly different result
from that in \cite{EG2}.\\
Let us see how the above lemma implies (i) and (ii).\\
\indent We start with the upper bound (i).\@ Take $a$, $\rho$, $T$
and $\delta$ positive.\@ Take $L>a$.\@ For $\tilde{a}$ in $(0,a]$,
we denote by
\begin{equation*}
A_{\tilde{a}}^{u_0}=\left\{v\in\xCn\left([0,T];\xHone\right):\
d_{\xCn\left([0,T];\xHone\right)}\left(v,K_T^{u_0}(\tilde{a})\right)
\geq\delta\right\}.
\end{equation*}
Note that we have $A_a^{u_0}\subset A_{\tilde{a}}^{u_0}$ and
$C_{\tilde{a}}\subset C_a$.\@ Take $\tilde{a}\in(0,a]$ and $f$ such
that $I_T^W(f)<\tilde{a}$.\\
\indent We shall now apply the Azencott lemma and choose $p=2$.\@ We
obtain $\epsilon_{\rho,f,\delta}$ and $\gamma_{\rho,f,\delta}$
positive such that for every $\epsilon\leq\epsilon_{\rho,f,\delta}$
and $u_0$ such that $\left\|u_0\right\|_{\xHone}\leq\rho$,
\begin{equation*}
\epsilon\log\mathbb{P}\left(\left\|u^{\epsilon,u_0}-\tilde{\mathbf{S}}(u_0,f)
\right\|_{X^{(T,p)}}\geq
\delta;\left\|\sqrt{\epsilon}W-f\right\|_{\xCn
\left([0,T];\xHn_{\xR}^s\right)}<\gamma_{\rho,f,\delta}\right)
\leq-L.
\end{equation*}
Let us denote by $O_{\rho,f,\delta}$ the set
$O_{\rho,f,\delta}=B_{\xCn\left([0,T];\xHn_{\xR}^s
\right)}(f,\gamma_{\rho,f,\delta})$.\@ The family
$\left(O_{\rho,f,\delta}\right)_{f\in C_a}$ is a covering by open
sets of the compact set $C_a$, thus there exists a finite
sub-covering of the form $\bigcup_{i=1}^NO_{\rho,f_i,\delta}.$ We
can now write
\begin{equation*}
\begin{array}{rl}
\mathbb{P}\left(u^{\epsilon,u_0}\in A_{\tilde{a}}^{u_0}\right)\leq&
\mathbb{P}\left(\left\{u^{\epsilon,u_0}\in
A_{\tilde{a}}^{u_0}\right\}\cap\left\{\sqrt{\epsilon}W\in\bigcup_{i=1}^N
O_{\rho,f_i,\delta}\right\}\right)\\
&+\mathbb{P}\left(\sqrt{\epsilon}W
\notin\bigcup_{i=1}^NO_{\rho,f_i,\delta}\right)\\
\leq &\sum_{i=1}^N\mathbb{P}\left(\left\{u^{\epsilon,u_0}\in
A_{\tilde{a}}^{u_0}\right\}\cap\left\{\sqrt{\epsilon}W\in
O_{\rho,f_i,\delta}\right\}\right)\\
&+\mathbb{P}\left(\sqrt{\epsilon}W\notin
C_a\right)\\
\leq
&\sum_{i=1}^N\mathbb{P}\left(\left\{\left\|u^{\epsilon,u_0}-\tilde{\mathbf{S}}(u_0,f)
\right\|_{X^{(T,p)}}\geq \delta\right\}\cap\left\{\sqrt{\epsilon}
W\in
O_{\rho,f_i,\delta}\right\}\right)\\
&+\exp\left(-a/\epsilon\right),
\end{array}
\end{equation*}
for $\epsilon\leq\epsilon_0$ for some $\epsilon_0$ positive.\@ We
used that
\begin{equation*}
d_{\xCn([0,T];\xHone)}\left(\tilde{\mathbf{S}}(u_0,f)
,A_{\tilde{a}}^{u_0}\right)\geq\delta,
\end{equation*}
which is a consequence of the definition of the sets $A_{\tilde{a}}^{u_0}$.\\
\indent As a consequence, for
$\epsilon\leq\epsilon_0\wedge\left(\min_{i=1,..,N}
\epsilon_{u_0,f_i}\right)$ we obtain for $u_0$ in $B_{\rho}^1$,
\begin{equation*}
\mathbb{P}\left(u^{\epsilon,u_0}\in A_{\tilde{a}}^{u_0}\right)\leq\
N\exp\left(-L/\epsilon\right)+\exp\left(-a/\epsilon\right),\end{equation*}
and for $\epsilon_1$ small enough, for every
$\epsilon\in(0,\epsilon_1)$,
\begin{equation*}
\epsilon\log\mathbb{P}\left(u^{\epsilon,u_0}\in
A_{\tilde{a}}^{u_0}\right)\leq\epsilon\log2+\left(\epsilon\log
N-L\right)\vee(-a).
\end{equation*}
If $\epsilon_1$ is also chosen such that
$\epsilon_1<\frac{\gamma}{\log(2)}\wedge\frac{L-a}{\log(N)}$ we
obtain
\begin{equation*}
\epsilon\log\mathbb{P}\left(u^{\epsilon,u_0}\in
A_{\tilde{a}}^{u_0}\right)\leq-\tilde{a}-\gamma,
\end{equation*}
which holds for every $u_0$ such that $\|u_0\|_{\xHone}\leq\rho$.\\
\indent We consider now the lower bound (ii).\@ Take $a$, $\rho$,
$T$ and $\delta$ positive.\@ The continuity of
$\tilde{\mathbf{S}}(u_0,\cdot)$, to be proved as in \cite{EG2},
along with the compactness of $C_a$ give that for $u_0$ such that
$\|u_0\|_{\xHone}\leq\rho$ and $w$ in $K_T^{u_0}(a)$, there exists
$f$ such that $w=\tilde{\mathbf{S}}(u_0,f)$ and
$I_T^{u_0}(w)=I_T^W(f)$.\@ Take $L>I^{u_0}(w)$.\@ Choose
$\epsilon_{\rho,f,\delta}$ positive and $O_{\rho,f,\delta}$, the
ball centered at $f$ of radius $\gamma_{\rho,f,\delta}$ defined as
previously, such that for every
$\epsilon\leq\epsilon_{\rho,f,\delta}$ and $u_0$ such that
$\left\|u_0\right\|_{\xHone}\leq\rho$,
\begin{equation*}
\epsilon\log\mathbb{P}\left(\left\|u^{\epsilon,u_0}-\tilde{\mathbf{S}}(u_0,f)
\right\|_{X^{(T,p)}}\geq
\delta;\left\|\sqrt{\epsilon}W-f\right\|_{\xCn
\left([0,T];\xHn_{\xR}^s\right)}<\gamma_{\rho,f,\delta}\right)
\leq-L.
\end{equation*}
We obtain
\begin{equation*}\begin{array}{rl}
\exp\left(-I_T^W(f)/\epsilon\right)&\leq\mathbb{P}\left(\sqrt{\epsilon}W\in
O_{\rho,f,\delta}\right)\\
&\leq\mathbb{P}\left(\left\{\left\|u^{\epsilon,u_0}-\tilde{\mathbf{S}}(u_0,f)
\right\|_{X^{(T,p)}}\geq
\delta\right\}\cap\left\{\sqrt{\epsilon}W\in
O_{\rho,f,\delta}\right\}\right)\\
&\quad+\mathbb{P}\left(\left\|u^{\epsilon,u_0}-\tilde{\mathbf{S}}(u_0,f)
\right\|_{X^{(T,p)}}<\delta\right).
\end{array}\end{equation*}
Thus, for $\epsilon\leq\epsilon_{\rho,f,\delta}$, for every $u_0$
such that $\|u_0\|_{\xHone}\leq\rho$,
\begin{equation*}
-I^{u_0}(w)\leq\epsilon\log 2+\left(\epsilon\log\mathbb{P}
\left(\left\|u^{\epsilon,u_0}-\tilde{\mathbf{S}}(u_0,f)
\right\|_{X^{(T,p)}}<\delta\right)\right)\vee(-L)
\end{equation*}
and for $\epsilon_1$ small enough and such that
$\epsilon_1\log(2)<\gamma$, for every $\epsilon$ positive such that
$\epsilon<\epsilon_1$, for every $u_0$ such that
$\|u_0\|_{\xHone}\leq\rho$,
\begin{equation*}
-I^{u_0}(w)-\gamma\leq\epsilon\log\mathbb{P}
\left(\left\|u^{\epsilon,u_0}-\tilde{\mathbf{S}}(u_0,f)
\right\|_{X^{(T,p)}}<\delta\right).
\end{equation*}
It ends the proof of (i) and (ii).\hfill $\square$

\end{document}